\documentclass[draft]{article}
\usepackage{textcomp}
\usepackage{amsmath}
\usepackage{amssymb}
\usepackage{amsthm}
\usepackage[final]{graphicx}

\DeclareMathOperator{\Int}{int}
\newtheorem{lemma}{Lemma}
\newtheorem{theorem}{Theorem}

\newsavebox{\picshift}
\newdimen\numheight
\sbox{\picshift}{1}
\numheight=\ht\picshift
\def\raisetonum#1{\sbox{\picshift}{#1}%
\raisebox{\numheight}{%
\raisebox{-\ht\picshift}{\usebox{\picshift}}%
}}

\newcommand{\N}{{\mathbb{N}}}
\newcommand{\Z}{{\mathbb{Z}}}
\newcommand{\R}{{\mathbb{R}}}

\newcommand{\eps}{\varepsilon}

\let\emptyset\varnothing

\DeclareMathOperator{\clos}{clos}
\DeclareMathOperator{\mes}{mes}

\newcommand\smallfour[4]{\begin{smallmatrix}#1&#2\\#3&#4\end{smallmatrix}}
\newcommand\psmallfour[4]{\bigl(\smallfour{#1}{#2}{#3}{#4}\bigr)}

\begin{document}

\begin{center}
\renewcommand{\thefootnote}{\fnsymbol{footnote}}
\textbf{\Large On the hyperbolic automorphisms of the 2-torus}\\
\textbf{\Large and their Markov partitions }\\[3mm]
\textbf{D.V. Anosov\footnote{Steklov Mathematical Institute,
Moscow, Russia, \texttt{anosov at mi.ras.ru}}, A.V.
Klimenko\footnote{Steklov Mathematical Institute, Moscow, Russia,
\texttt{klimenko05 at mail.ru}}, G. Kolutsky\footnote{Lomonosov
Moscow State University, Moscow, Russia, \texttt{kolutsky at
mccme.ru}}}
\end{center}

\setcounter{footnote}{0}
\renewcommand{\thefootnote}{\arabic{footnote}}

\section*{Key words and phrases}

hyperbolicity, symbolic dynamics, Anosov maps, Markov partitions

\section*{Abstract}

   An (algebraic) automorphism of the 2-torus is defined in a standard way by a
matrix with determinant 1 or $-1$ and with integer coefficients.
An automorphism is hyperbolic, if the eigenvalues of this matrix
are reals with absolute value $> 1$ for one eigenvalue (and $<1$
for another). Iterations of such automorphism $A$ constitute a
dynamical system (DS) with discrete time~--- phase points do not
move continuously as it is for the DS described by differential
equations, but jump from one place to another; the moving phase
point which originally (at the zero moment of time) occupied the
position $x$ moves to $A^n x$ during the time $n$. Hyperbolicity
implies that although formally this DS is deterministic, actually
the behavior of its trajectories resembles, in a sense, behaviour
of some random (stochastic) process. Markov partitions is the best
method to establish this analogy which is even a kind of
isomorphism.

     This text is based on the talk
the first author gave in Germany, but the text is more detailed.
It consists of four parts.\footnote{ In the lecture he was
restricted in time. However, here we also omit some details. Still
we think that the mainstream is more or less clear and that a
competent mathematician can easily elaborate the omitted details
belonging to the mainstream.} In the first part we explain how the
deterministic DS can be isomorphic to a random process on an
example (the circle expanding map) which is more simple. In the
second part we dwell on the classification of hyperbolic toric
automorphisms. In the third part we define the notion of Markov
partitions and explain how they can be used and how one can
construct a simplest Markov partition (perhaps some details of the
construction can be somewhat new). Finally, in the fourth part we
describe a kind of classification of these simplest Markov
partitions (this is new).

Parts 2, 3 and 4 are based on the work of A.V.\,Klimenko and G.
Kolutsky who are Ph.D.~students of D.V.\,Anosov. Besides him, in
the beginning of their work their inofficial scientific advisor
was A.Yu.\,Zhirov. Part 2 is an exposition of results which seems
to be known in the number theory; the version presented here was
elaborated by G. Kolutsky. Parts 3 and 4 is mainly due to
Klimenko; the idea of using results and notions from the Part 2
for the goals of Part 4 was a result of his discussion of the
matter with Kolutsky; also, they examined several first examples
together.

The first author thanks the Humboldt Foundation, which supported
his visit to Germany in 2007, and the Max Planck Institute for
Mathematics and the Ulm University for their hospitality. Our work
exposed here was supported by the grants No. 05-01-01004-a and
08-01-00342-a of the Russian Foundation for Basic Research, by the
grants No. NSh-6849.2006.1 and NSh-3038.2008.1 of the President of
Russia for support of leading scientific schools and by the
Program ``Nonlinear dynamics'' of the Russian Academy of Sciences.

\section{Introduction}

   Two big parts of the theory of dynamical systems can be characterized
as dealing with motions of ``regular'' and ``stochastic, quasi-random,
chaotic'' character.  Simplest examples of regular motions (and those which
are, informally, ``the most regular'') are periodic or quasiperiodic motions.
(Thus considering of regular motions is as old as the science itself ---
some regularity of planets' motion was known and exploited by Babylonians,
and in the more advanced Ptolemeus' system these motions were essentially
described by trigonometric polynomials.) Examples of ``chaotic'' motions
are much more new. As far as we know, the first example of such kind was
pointed out by J.~Hadamard about 1900. A couple of decades earlier H.Poincar\'e
discovered the so-called ``homoclinic points'' which now serve as practically
the main ``source'' of ``chaoticity''; however, Poincar\'e himself spoke only
that the ``phase portrait'' (i.~e.\ the qualitative picture of trajectories'
behaviour in the phase space) near such points is extremely complicated.
A couple of decades after Hadamard E.Borel encountered a much simpler
example of the ``chaoticity'' where it is easy to
understand the ``moving strings'' of this phenomenon. We shall begin with
a description of his example. About 100 years later it remains the simplest
manifestation of the fact that a dynamical system (which, by definition, is
deterministic) can somehow resemble a stochastic process (in fact, even be,
in a reasonable sense, isomorphic to such process).

In this example the phase space is the circle $\mathbb{S}^1 = \mathbb{R}/\mathbb{Z}$.
We shall often speak that $\mathbb{R}$ projects onto $\mathbb{S}^1$ by the projection $p$.
We can consider the usual coordinate $x$ in $\mathbb{R}$ as a ``cyclic coordinate''
on $\mathbb{S}^1$. In its terms we define the map
\begin{equation}\label{DoubleMap}
f\colon\mathbb{S}^1 \to \mathbb{S}^1, \qquad f(x) = 2x.
\end{equation}
More formally, we begin with the map
\begin{equation*}
g\colon\mathbb{R} \to \mathbb{R}, \qquad x  \mapsto 2x
\end{equation*}
and project it onto $\mathbb{S}^1$ (so $p(x) \mapsto p(2x)$; we use the fact
that points $2x$ and $2(x + n)$ ($n$ is an integer) project to the same point
of $\mathbb{S}^1$. More formally, we use that $g(\mathbb{Z})\subset\mathbb{Z}$
so that $g$ maps the class $x+ \mathbb{Z}$ to the class $2x+ \mathbb{Z}$.)
Pictorially, considering $\mathbb{S}^1$  as made from rubber, we stretch it
to double its length and then cover the original $\mathbb{S}^1$ by this
expanded circle (so each point of the initial circle is covered by two points
of the expanded one).\footnote{This $f$ is an example of the so-called
``expanding diffeomorphism'' of $\mathbb{S}^1$. We shall not need to define
this class of maps, as we shall deal with $f$ only. But on the ``conversational
level'' it is clear that $f$ deserves to be called ``expanding''.}

Our dynamical system consists of iterations $\{f^n\}$ of $f$, so that any of
its trajectories is a sequence $\{f^n(x), \ n \in \mathbb{Z}_+\}$
(here, in Bourbaki's style, $+$ is used to deceive a spy; actually
$\mathbb{Z}_+ = \{0,1,2,\ldots\}$). Thus it is a system with discrete time
($n$ plays the role of time --- during time $n$ the moving phase
point ``jumps'' from the original position $x$ into the position
$f^n(x)$).

Remark: One can inquire whether it is possible to construct a
system with continuous time exhibiting ``chaotic'' properties
analogous to those we are going to discuss for our $\{f^n\}$; and
whether there exist dynamical systems with chaotic behavior of
their trajectories among those systems of the most classical
character --- those described by phase velocity vector fields $v$
on a smooth phase manifold~$M$ (the moving phase point moves
accordingly to the differential equation $\dot x = f(x)$ which in
terms of local coordinates looks as a ``habitual'' system of
autonomous differential equations). The answer is positive.
Essentially first examples of such kind were found in the process
of improving Hadamard's results. But for all such systems the
phase space is unavoidably of dimension not less than $3$ and they
are much more complicated than Borel's example.

Another question preceding discussion of any concrete properties
of Borel's example is the following. In this example the map $f$
is irreversible; so we can speak about the future motion of the
moving phase point (it occupies position $x$, then $f(x)$, them
$f^2(x)$, and so on), but we can't speak about its position for
negative time $n$. Is it possible to construct ``reversible
chaotic'' examples? Basically the positive answer to the previous
question indicates that this is possible (in ``classical''
dynamical systems the time is reversible), so that the question
can be only for dimension of the phase space less than 3. This can
be achieved if we pass from the continuous time to a discrete one.
Actually the main content of this paper will be related to the
simplest example of such kind. Reversibility is gained at the
price of increasing the phase space dimension~--- 2 instead of 1;
namely, we shall deal with a smooth automorphism of the 2-torus.
But we begin with Borel's example, as it is more simple.

From now on till the end of this part $f$ means Borel's $f$ defined
by~\eqref{DoubleMap}. If we knew $x$
precisely, we could compute its trajectory $\{f^n(x)\}$. But assume that
we know the phase point we have to deal with only approximately, although
with a good approximation.
So instead of the ``true'' trajectory $\{f^n(x)\}$ (or $\{2^n x\}$ in terms of the
cyclic coordinates) we compute the trajectory $\{f^n(y)\} = \{2^ny\}$ with some
$y$ at the small distance $\delta$ from $x$.
The distance between $f^n(x)$ and $f^n(y)$ is $2^n\delta$. For
several first numbers $n$ the error is small, but it rapidly increases with $n$.
Without entering into refinements of the terminology, this can be called instability,
and even a strong one --- roughly speaking, this kind of instability means that
two phase points which originally were close to each other can rapidly diverge
under the action of the iterations $f^n$. (More technically, such type
of instability
is called exponential, uniform and complete; we shall not dwell on this.)
If $\delta$ is of the order $10^{-8}$ (the size of atom in centimeters), for
$n = 30$ the error will be of order $10$, i.~e.\ of the macroscopic order~--- of
the same size as the laboratory equipment or (returning to our example) as our
circle $\mathbb{S}^1$ (formally, even more than it). Then all what we can say
is that the moving phase point $f^n(x)$ is somewhere on the circle --- a
trivial conclusion which can be made without any measurements and calculations.

Besides this ``growth of uncertainty'' which comes to attention
when we compare the  behaviour of two different trajectories
$f^n(x)$ and $f^n(y)$ (with $x \approx y$), behavior of the most
part of individual trajectories $f^n(x)$ also demonstrates such
features which make it reasonable to characterize their behavior
as a chaotic one. We shall see this later.

 About 1910 Poincar\'e wrote that in such situation instead of the more or less
exact computing the ``individual'' trajectory (which is practically impossible)
one can try to make some statistical statements concerning some features of
behaviour of a ``majority'' of trajectories or of the ``typical''
trajectories. Instability, in his opinion, was the source (which can be
a hidden source) of the probability.

   We suspect that besides Poincar\'e some physicists also shared this point
of view at that time (very end of XIX~--- beginning of XX century). But, in
any case, he expressed it quite distinctively and illustrated it on some
mathematical example. We shall not dwell on it because the later Borel's example
provides a better illustration which at the same time is more close to the
goal of this paper. (In Poincar\'e's example individual trajectories were not
chaotic and the distance between $f^n(x)$ and $f^n(y)$
was growing more slowly than in Borel's case.)

    Now we know that besides instability there exists at least one source
of the random behavior, that is, quantum effects. But this does not abolish
those effects which are due to the instability and so emerge even in
the classical situation.

   Actually Borel spoke not about the circle map $f$, but about the interval map
\begin{equation*}
[0,1) \to [0,1), \qquad x \mapsto \{2x\} \quad \text{($\{\,\cdot\,\}$ means the
fractional part).}
\end{equation*}
This map has a disadvantage of being discontinuous at the point $x
=1/2$. For the reason to be explained below this discontinuity did
not trouble Borel. However, we see that we can easily get rid of
it --- just replacing $[0,1)$ by $\mathbb{S}^1$.

   In the original Borel's version it is especially clear that the map $f$ is quite
lucidly described in terms of the expansion of $x$ into infinite binary fraction,
If, in these terms,
\begin{equation*}
x = 0{,}a_1a_2a_3\ldots \qquad \text{with all $a_i$ being 0 or 1,}
\end{equation*}
which means that
\begin{equation*}
x = \frac{a_1}{2} + \frac{a_2}{2^2} + \frac{a_3}{2^3} + \ldots,
\end{equation*}
then $f(x) = 0{,}a_2a_3a_4\ldots$. The comma separating the integer part of
the binary fraction from its fractional part is moved one step to the right
and all that becomes to the left of the shifted comma is replaced by zero.
One can also say that the comma's position is fixed, but the infinite sequence
$a_1a_2a_3\ldots$ shifts by the one step to the left and the coefficient
$a_1$ (appeared to be to the left of the comma) is discarded (i.~e.\ replaced
by 0). The binary expansion of $x$ is not unique
for binary-rational $x$ (e.g. for those of the form $x = \frac{\text{integer}}{2^n}$).
But it is harmless, because if two binary expansions represent the same $x$,
the shifted binary expansions represent the same $f(x)$.

   In terms of the circle $\mathbb{S}^1$ one can interpret the binary
expansions as follows. Points  $p(\frac{i}{2^n}) \quad (i = 0,\ldots,2^n - 1)$
divide $\mathbb{S}^1$
into $2^n$ arcs. $(i + 1)$-th arc consists of points $p(x)$ obtained when $x$
increases from $\frac{i}{2^n}$ to $\frac{i + 1}{2^n}$; i.e., this arc is
$p\left(\left[ \frac{i}{2^n},\frac{i + 1}{2^n}\right]\right)$. Let us denote this arcs
as follows. If $b_k\ldots b_1b_0$ is the binary representation of $i$,
we define $b_{k+1}=b_{k+2}=\dots=b_{n-1}=0$ and then associate with
each $i=0,1,\dots, 2^n-1$ the sequence $b_{n-1},\dots, b_0$.
E.~g., binary representation for $i = 3$ is 11, and if $n = 4$, we associate
with $3$ the finite sequence 0011.) Having in mind this correspondence between
numbers $i$ and sequences
$b_{n - 1}\ldots b_0$, denote
\begin{equation*}
p\left(\left[ \frac{i}{2^n},\frac{i + 1}{2^n}\right]\right)
= C_{b_{n - 1}\ldots b_0}.
\end{equation*}
Then\footnote{Here and below $*$ denotes an arbitrary digit.}
\begin{equation*}
p(x) \in C_{b_{n - 1}\ldots b_0}\quad\text{ if and only if }
x = 0{,}b_{n - 1}\ldots b_0 *\ldots *\ldots.
\end{equation*}
A point with binary rational cyclic coordinate has two binary expansions ---
say,
\begin{equation}\label{TwoBinExp}
0{,}a_1\ldots a_k 0 1\ldots 1 \ldots \quad\text{and}\quad
0{,}a_1 \ldots a_k 1 0 \ldots 0 \ldots.
\end{equation}
If $k \ge n$, first $n$ coefficients of these expansion are the same, and so
for both expansions our receipt says that $p(x) \in C_{a_1 \ldots a_n}$. If
$k <n$, the point $p(x)$ is the endpoint of two adjacent arcs $C_{c_1\ldots c_n}$,
and their labels $c_1\dots c_n$ will be first $n$ digits of one or another
binary expansion~\eqref{TwoBinExp}.

This geometric characterization of the binary expansion of $x$ is,
so to say, a ``static'' one. But it is easy to pass to a
``dynamical'' characterization of this expansion:
\begin{equation*}\arraycolsep=0pt
\begin{array}{rll}
x={}&0{,}a_1 *\ldots * \ldots&\mbox{\ \ if and only if } p(x) \in C_{a_1},\\
x={}&0{,}a_1 a_2 *\ldots * \ldots&\mbox{\ \ if and only if } p(x) \in C_{a_1},
f(p(x)) \in C_{a_2},\\
\multicolumn{3}{l}{\text{\indent(recall that $f(p(x)) = 0{,}a_2 * \ldots * \ldots $);}}\\
\multicolumn{3}{l}{\dotfill\quad}\\
x={}&0{,}a_1\ldots a_n *\ldots * \ldots&\mbox{\ \ if and only if }
p(x) \in C_{a_1},\ldots, f^{n - 1}(p(x)) \in C_{a_n},\\
\multicolumn{3}{l}{\dotfill\hskip -0.2em}\\
\end{array}
\end{equation*}

   Of course it is only the sequence $\{a_n\}$ that is important,
not the zero and comma standing before them. Slightly modifying what was said
earlier (and deviating from literally following Borel), we can adopt the
following agreements. Instead of numbers $x \in [0,1)$
we shall begin with (singly-) infinite sequences $(a_0,\ldots,a_n,\ldots)$ of
numbers (or symbols) $a_i \in \{0,1\}$ (now we start numbering them
from 0; advantage of this is that now $a_n$ is the number of the semicircle
$C_i$ containing $f^n(x)$). Denote by $\Omega$ the space of all these
sequences (i.e., $\Omega = \{0,1\}^{\mathbb{Z}_+}$). Word ``space''
hints that  $\Omega$ will not be merely a set, but that it will be
endowed with some structure. There will be two structures on $\Omega$:
topology and measure.

   As regards to topology, we take the discrete topology (each point is an open set) in each
multiplier $\{0,1\}$ of the infinite product $\{0,1\}^{\mathbb{Z}_+}$ and then endow
this product by the Tikhonov product topology. According to Tikhonov theorem,
$\Omega$ is compact as a product of compact spaces. In this case the topology
on $\Omega$ is induced by some metric, e.g. one can take
\begin{equation*}
\rho (x,y) = \sum_n \frac{ d(x_n,y_n)}{2^{n+1}} \quad
\mbox{for } x = (x_0,x_1,\ldots), \ y = (y_0,y_1,\ldots),
\end{equation*}
where $d(a,b) = 0$ for $a = b$ and $d(a,b) = 1$ for $a \neq b$. Using this
metric, one can easily prove compactness of $\Omega$ without referring to the
general theorem.

   Subset $A \subset \Omega$ is called a cylindric set if it consists of all
sequences $x$ such that some prescribed coordinates $x_{i_1},\ldots,x_{i_n}$
of $x$ are given numbers $a_{i_1},\ldots,a_{i_n}$, while other coordinates are
arbitrary. Cylindric sets are open in the topology used; moreover, they
constitute a base for this topology. They are also closed --- existence of so many
open-closed sets means that $\Omega$ is zero-dimensional.

 As we've started to speak about products, we shall sometimes call the $n$-th
element $x_n$ of the sequence $x = (x_0,x_1,\ldots)$ its $n$-th coordinate
(once more, they are numbered beginning from the 0-th coordinate).

   Binary expansions were binary expansions of the cyclic coordinates of the
points of $\mathbb{S}^1$. In our new language we introduce the map
\begin{equation}\label{pitoS1}
\pi: \Omega \to \mathbb{S}^1 \qquad \pi(x) =
p\left(\sum_n \frac{x_n}{2^{n + 1}}\right).
\end{equation}
It is a continuous map. There exist a countable set of points having two
preimages, but for the ``vast majority'' of points there is only one preimage.
Multiplying cyclic coordinates by 2 is now replaced by the ``one-side Bernoulli
shift'' $\sigma$ moving the whole sequence to one step left and omitting its
first symbol; that is,
$$ \mbox{for } x = (x_0,x_1,\ldots) \quad \sigma(x) = (y_1,y_2,\ldots),
\quad \mbox{where } y_n = x_{n + 1} \mbox{ for all } n \in \mathbb{Z}_+.$$
It is clear that $\pi \circ \sigma = f \circ \pi$. In this sense one can say
that our construction provides a ``symbolic model'' for our original map
$f:\mathbb{S}^1 \to \mathbb{S}^1$.

   Point $x$ and its trajectory $\{f^n(x)\}$ are ``coded'' by a sequence
$(a_0,a_1,\ldots)$ (once more: $n$-th element of this sequence is such number
that $f^n(x) \in C_{a_n}$). This sequence could be called ``a journey diary
of $x$''. Yu.S.Il'yashenko uses the more impressive name ``a fate of $x$''.
Below we often call this sequence simply ``a code of $x$''.

     This trick --- ``diary'', ``fate'', ``coding'' ---
is by no means restricted by our example. If some set $X$ is decomposed
into nonintersecting sets
\begin{equation}\label{decomp}
X = X_1 \cup \ldots \cup X_k, \qquad X_i \cap X_j = \emptyset
\mbox{ for } i \neq j,
\end{equation}
then for any map $f: X \to X$ we can introduce ``a journey diary'' of a
point $x \in X$ (with respect to the decomposition (\ref{decomp})):
this ``diary'' is an infinite sequence $(a_n; \ n \in \mathbb{Z}_+)$ such that
$f^n(x) \in X_{a_n}$. Of course, the decomposition (\ref{decomp}) must be
somehow adjusted to the structures which are specific for example or a
class of examples we are going to consider (and which are somehow respected
by $f$). Besides this general demand, a special choice of the decomposition
used may take into account more specific properties of $f$. Also, in our
case this general approach is slightly modified.
Essentially we are using the partition $\mathbb{S}^1 = C_0 \cup C_1$ which is
not a decomposition in the strict sense: $C_0 \cup C_1 \neq \emptyset$.
As a result, the encoding the point $x$ by sequence $(a_n)$ does not always
supply us with a single valued function $x \mapsto (a_n)$:
some points of $\mathbb{S}^1$ (those with binary-rational cyclic coordinates)
have several (two) ``journey diaries''. This would not happen if we took
$C_0 = p\left(\left[0,\frac{1}{2}\right)\right), \ C_1 =
p\left(\left[0,\frac{1}{2}\right)\right)$. On the language of the binary
expansions, this would mean that we rule out expansions of the form
$0,a_1\ldots a_k 11\ldots 1 \ldots$, i.e. those to be periodic after
some place with the period\footnote{Here and later we shall often use the
word ``period'' as denoting the periodic part of the infinite sequence,
not merely the length of this part.} consisting of one digit 1. However,
practically one uses such binary expansions and we shall also use the closed
arcs $C_i$.

     Our ``journey diary'' can be described in accordance to a general remark
above in terms of dynamics and partition $\mathbb{S}^1 = C_0 \cup C_1$, without
appealing to binary expansions:
\begin{equation}\label{coding}
 x \mapsto (a_n) \quad \mbox{if and only if } f^n(x) \in C_{a_n}
\quad \mbox{for all } n \in \mathbb{Z}_+.
\end{equation}
This makes evident that if $x \mapsto a = (a_0,a_1,a_2,\ldots)$, then $f(x) \mapsto
(a_1,a_2,a_3,\ldots)$. But essentially we have also used the binary expansions
in the definition of the map (\ref{pitoS1}) inverse to the (multi-valued)
coding $x \mapsto (a_n)$ (which makes it evident that any sequence $(a_n)$
codes some $x$). Here it is also easy to get rid of them. (\ref{coding}) is
equivalent to
$\pi((a_n)) \in \bigcap_{n = 0}^{\infty} f^{-n}(C_{a_n})$, i.e.
\begin{equation}\label{defpi}
 \mbox{for all } N \in \mathbb{Z}_+ \quad \pi((a_n)) \in \bigcap_{n = 0}^N
f^{-n}(C_{a_n}).
\end{equation}
Define $F_N = \bigcap_{n = 0}^N f^{-n}(C_{a_n})$. Clearly $F_0 \supset F_1 \supset
\ldots \supset F_N \supset \ldots$. It turns out that
\begin{equation}\label{arcFN}
F_N \quad \mbox{is a closed arc of the length} \quad \frac{1}{2^{N + 1}}.
\end{equation}
This implies existence and uniqueness of the point common to all $F_N$.
This implies also the continuity of $\pi$. Indeed, if $\rho((a_n),(b_n))$
is small, which implies that $a_n = b_n$ for all $n = 0,1,\ldots,N$ with
some big $N$, then both $\pi((a_n))$ and $\pi((b_n))$ lie within the same
arc $F_N$ of the small length $\frac{1}{2^{N + 1}}$.

     As regards to (\ref{arcFN}), it can be proved as follows. Clearly
$f^{-n}(C_0)$ and $f^{-n}(C_1)$ are disjoint unions of $2^N$ closed arcs
of the view $\left[\frac{i}{2^{n + 1}},\frac{i + 1}{2^{n + 1}}\right]$
with some $i \in \{0,1,\ldots,2^{n + 1} - 1\}$, $i$ being even for arcs
from $f^{-n}(C_0)$ and odd for arcs from $f^{-n}(C_1)$ ($f^n$ maps
homeomorphically any such arc with an even $i$ onto $C_0$ and with an odd
$i$ --- onto $C_1$). Any arc
$\left[\frac{i}{2^{n + 1}},\frac{i + 1}{2^{n + 1}}\right]$ consists of two
arcs of the form
\begin{equation}\label{2subarcs}
 \left[\frac{2j}{2^{n + 2}},\frac{2j + 1}{2^{n + 2}}\right], \quad
     \left[\frac{2j + 1}{2^{n + 2}},\frac{2j + 2}{2^{n + 2}}\right].
\end{equation}
Thus if we already know that $F_N$ is an arc of the type described (which
is trivial for $N = 0$), then passing to $F_{N + 1}$ means that we pass to
one of the arcs (\ref{2subarcs}) (to the first arc if $a_{N + 1} = 0$
and to the second arc if $a_{N + 1} = 1$).

   Our map $f$ is very simple, and at the first glance it is not clear whether our
symbolic model is useful for any purpose. We shall see that it is.

    It turns out that one can introduce a measure $\mu$ on $\Omega$
such that $\mu(A) = \frac{1}{2^n}$ for any cylindric $A$ defined by fixing
$n$ coordinates. (Of course dealing with the topological space we consider only
measures which are in a sense compatible with topology. In our case when the
space is a metrizable compact set this means simply that all Borel sets are
measurable.) Existence of such measure is a simple case of some general theorems
of the measure theory  and/or of the probability theory, but in this case
argumentation can be much more easy. Consider first the cylindric sets of
the following special character: they are defined by fixing first $n$
coordinates of their points; i.e. we speak about the sets
\begin{equation*}
B_{a_0,\ldots,a_{n-1}} = \{x = (x_0,x_1,\ldots); \ x_0=a_0, \ldots,
x_{n-1}=a_{n-1}\}.
\end{equation*}
This set is mapped under $\pi$ on the arc $C_{a_0,\ldots,a_{n - 1}}$.
The length of this arc is equal to $1/2^n$ which is just what we want to
be the measure of $B_{a_0,\ldots,a_{n - 1}}$. Going further, we observe
that any cylindric set $A$
is a finite union of the sets $B_{a_0,\ldots,a_{n - 1}}$ and $\pi$ maps such
union onto a finite system of arcs considered. It is easy to check that
the total length of these arcs is just what we want to be $\mu(A)$. And
this gives us an idea how to define $\mu$: we simply define it as the preimage
of the standard Lebesgue measure (denoted by $\mes$) on $\mathbb{S}^1$ (or, if you
prefer, on $[0,1)$ --- the Lebesgue measure does not feel the difference
between them which is due to just one point) under the map $\pi$.
Although $\pi$ is not a bijection, the violation of bijectivity is negligible
from the measure-theoretic point of view. So $\pi$ is an isomorphism of the
measure spaces $(\Omega,\mu)$ and $(\mathbb{S}^1,\mbox{mes})$.

   An important property of this measure is that for any measurable set
$A \subset \Omega$ its preimage $\sigma^{-1}(A)$ is also measurable (thus
$\sigma$ is measurable) and
\begin{equation}\label{pres}
\mu(\sigma^{-1}(A))=\mu(A).
\end{equation}
In such cases one says that the measure $\mu$ is invariant with
respect to $\sigma$. (Literally this expression would mean that
$\mu(\sigma (A)) = \mu(A)$. But this is wrong. When dealing with
any noninvertible map $\sigma$, one always understands
preservation of measure as the measurability of this map plus the
property (\ref{pres}).)

{\footnotesize Basic fact here is that these two properties
(measurability of $\sigma^{-1}(A)$ and (\ref{pres}))
are true for cylindric $A$. Let $A$ be described by fixing coordinates
$x_{i_1},\ldots,x_{i_n}$ of its points $x$ (so $\mu(A) = \frac{1}{2^n}$).
Preimage $\sigma^{-1}(x)$ consists of two points $y$ and $z$. Both have
the same coordinates which number is $i>0$ --- namely,
$y_i = z_i = x_{{i - 1}}$ (indeed, after the shift of $y$ and $z$
towards one step to the left one must get $x_{i - 1}$ on the
$(i - 1)$-st place), while $y_0 = 0$ and $z_0 = 1$ (thus no restrictions are
imposed on the zero's coordinate of the points of $\sigma^{-1}(A)$ ---
it can be 0 or 1 and this has no influence on other coordinates). It
follows that  $\sigma^{-1}(A)$ is the cylindric set such that restrictions
on the coordinates are imposed on the coordinates
$x_{i_1 + 1},\ldots, x_{i_n + 1}$. This is $n$ coordinates and so
$\mu(\sigma^{-1}(A))= \frac{1}{2^n} = \mu(A).$

   After this one can use more or less standard arguments from the measure
theory. We shall repeat them making simplifications due to specific features
of our case. Let $A$ be the finite union of cylindric sets $A_1,\ldots,A_n$.
Then $\sigma^{-1}(A)$ is a finite union of their preimages $\sigma^{-1}(A)$
which are also cylindric sets and thus measurable. This proves the measurability
of $\sigma^{-1}(A)$. Comparison of its measure with the measure of original
$A$ needs more considerations. Each $A_i$ is described by fixing a finite number
of coordinates --- say, fixing coordinates $x_j$ with $j \in J_i$ where $J_i$
is some finite set of nonnegative integers. Let $N =
\max (J_1 \cup \ldots \cup J_n)$. Any $A_i$ can be presented as a finite
union of some sets of the form $B_{a_0,\ldots,a_N}$. (Say, let the
restrictions describing $A_1$ be $x_1 = 0, x_2 = 1$ and the restrictions
describing $A_2$ be $x_0 = 1$ and $x_4 = 0$. Then $J_1 \cup J_2 =
\{0,1,2,4\}$ and $N = 4$.
We have
\begin{align*}
A_1={}&B_{00100} \cup B_{00101} \cup B_{00110}\cup B_{00111} \cup
B_{10100} \cup B_{10101} \cup B_{10110}\cup B_{10111},\\
A_2={}&\mbox{union of 8 sets $B_{1,a_1,a_2,a_3,0}$ for all $(a_1,a_2,a_3)\in\{0,1\}^3$.}
\end{align*}
Finite union of $A_i$ is also a finite union of some $B_{a_0,\ldots,a_N}$.
As these $B_{\ldots}$ do not intersect each other and
$\mu(\sigma^{-1}(B_{a_0,\ldots,a_N})) = \mu(B_{a_0,\ldots,a_N})$, it follows
that $\mu(\sigma^{-1}(A)) = \mu(A)$.

   Now any open set $U$  can be represented as a union of increasing sequence
$$U_1 \subset U_2 \subset \ldots \subset U_n \subset \ldots$$
of the sets each of which is a finite union of cylindric sets. (So $
\mu(U) = \lim\limits_{n\to\infty} \mu(U_n)$.) Then $\sigma^{-1}$ is the union
of increasing sequence
$$\sigma^{-1}(U_1) \subset \sigma^{-1}(U_2) \subset \ldots
\subset \sigma^{-1}(U_n) \subset \ldots.$$
Each $\sigma^{-1}(U_n)$ is measurable
(thus the union $\sigma^{-1}(A)$ of these sets is also measurable and
$\mu(\sigma^{-1}(A)) = \lim\limits_{n\to\infty}\mu( \sigma^{-1}(U_n))$) and has the same measure
as $U_n$. It follows that $\mu(\sigma^{-1}(A)) = \mu(A)$.

   Next step is to consider closed $A$. As $\sigma^{-1}(A) =
\Omega \setminus \sigma^{-1}(\Omega \setminus A)$, it is easy to see that
$\sigma^{-1}(A)$ is measurable and its measure equals to $\mu(A)$.

   Finally consider arbitrary measurable $A$. For any $\varepsilon > 0$ there
exist a closed set $C$ and an open set $U$ such that $C \subset A \subset U$
and $\mu(U) - \mu(C) < \varepsilon$ (in particular, $|\mu(U) - \mu(A)| <
\varepsilon$). Then $\sigma^{-1}(C) \subset \sigma^{-1}(A) \subset
\sigma^{-1}(U)$, the first set is closed, the last set is open and the
difference of their measures is the same as for original $U,C$, i.e. it
is less than $\varepsilon$. The fact that $\sigma^{-1}(A)$ contains some
measurable set and is contained in some open set and the measures of these
sets can be made arbitrarily close to each other, implies that $\sigma^{-1}(A)$
is measurable. It follows also that $|\mu(\sigma^{-1}(A)) - \mu(\sigma^{-1}(U)|
< \varepsilon$. And as $\mu(\sigma^{-1}(U)) = \mu(U)$, we see that
$|\mu(\sigma^{-1}(A) - \mu(A)| < 2\varepsilon$. As $\varepsilon$ is arbitrary,
we conclude that $\mu(\sigma^{-1}(A)) = \mu(A)$.}

   Now it is time to explain what was discovered by Borel (not the description
of the multiplication by 2 in terms of binary expansions, of
course). Borel observed that the dynamical system
$(\Omega,\sigma,\mu)$ \footnote{As we have already said, actually
he spoke of $([0,1), x\mapsto\{2x\}, \mes)$, but this difference
is not important from the point of view of his goal.} describes
the classical object of the probability theory~--- a sequence of
independent trials consisting in flipping of a coin. This
discovery was important for the development of the treatment of
probability theory foundations on the base of measure
theory\footnote{Borel's work was also influential in other
respects (some hint on this will be given below), but at the
moment we dwell only on one side of it which is close to our main
topic.}. In full generality this treatment was elaborated by
A.N.Kolmogorov in 1930s and became standard. Having this treatment
in mind, we can consider $(\Omega,\sigma,\mu)$ as an early
manifestation of this treatment applied to the coin flippings.

   We shall use three basic notions: a random event, probability and
independence. Essentially they cannot be defined in terms of
notions from other parts of the science. They can be only
illustrated on examples on semi-intuitive level. But the mutual
relations of these notions can be described completely using other
mathematical notions. Essentially this is the usual situation with
basic notions in any part of mathematics\footnote{Euclidus' claim
that ``a point is what has no parts'' so often criticized as
``naive, obscure and having no real content'' is merely a naive
way to say that in Euclidean geometry we deal with some sets
(3-dimensional Euclidean space and its subsets) endowed with some
structure described by the axioms and that points are just
elements of these sets. As those, they really have no parts,
Hilbert space $H$ can well be some class of functions and
functions themselves are rather complicated things; but as a point
of $H$ each function is considered as something what is
``primitive, elementary, without intrinsic structure''.}.

  First consider finite sequences of independent coin flippings. Say, let
us flip a coin three times. An example of the random event: we
have got 0 after the first flip, 1 after the second flip, and 0
after the third one. This can be denoted by the finite sequence
(0,1,0). This is an example of what is called an elementary event.
In our case the elementary event describes the result of a
flipping repeated three times. So there are eight elementary
events described by 8 binary sequences $(a_1,a_2,a_3)$ with all
$a_i = 0$ or 1. We can even adopt a formal point of view
considering these sequences themselves as elementary events. Their
collection $\{0,1\}^3$ is what is called the space of elementary
events. An example of a non-elementary event $A$: the sum of the
numbers associated with three flips is odd. This happens if and
only if the results of three subsequent coin flips are $(0,0,1),
(0,1,0), (1,0,0), (1,1,1)$. Thus we can consider an event as a
subset of the space of elementary events. An event $B$ consisting
in 0 being the result of the first flip and the sum of the numbers
associated with 3 flips being odd is a subset of the previous $A$
consisting of $(0,0,1)$ and $(0,1,0)$. Going further, we say that
any result of a single flip of the coin appears with the
probability $\frac{1}{2}$. (This is practically interpreted that
if we flip the coin many times or if we flip many coins
simultaneously, approximately half of these trials will have the
result 0. Once more: from the point of view described this
statement is not the definition of the probability, but merely a
kind of intuitive explanation, or illustration, of this basic
notion.) It is because the coin is assumed to be ``fair'', i.~e.\
symmetric with respect to both its sides. Independence of the
subsequent flips of the coin manifests itself in the fact that
probability of any elementary event $(a_0,a_1,a_2)$ is
$\frac{1}{2^3}$.

   We do not know whether there exist ``false'' coins such that the probabilities
of 0 and 1 are considerably different from
$\frac{1}{2}$.\footnote{There are similar procedures with
probability different from $1/2$. For example, spinning of a
newly-minted U.S. penny on a smooth table tends to show less
``heads'' than ``tails'' (as Lincoln's head overweighs another
side). For some manners of spinning the probability of ``head''
can be as small as $0{,}1$.} But there certainly exist loaded
dices. According to the literature, they are even of some
practical importance. If the dice is ``fair'', i.e. symmetric with
respect to its faces and made from homogeneous material, then  the
probability of any of its faces to be shown after throwing of the
dice is $\frac{1}{6}$. For loaded dice they are some numbers
$p_1,p_2,p_3,p_4,p_4,p_6$ such that all $p_i \geq 0$ and $\sum p_i
= 1$. Assuming that we deal with a nonsymmetric coin, there is a
probability $p_0$ that the result of a flip of the coin will be 0
and a probability $p_1$ that this result will be 1. Numbers $p_i$
are $\geq 0$ and their sum $p_0 + p_1 = 1$. In such case an
elementary event $(a_1,a_2,\ldots,a_n)$ has the probability
$p_{a_1}p_{a_2}\ldots p_{a_n}$.

   Be the coin fair or not, after we defined the probabilities of
elementary events, probability of any event $A$ is just the sum of probabilities
of its elements (of the elementary events belonging to $A$). So we get some
structure on the space of elementary events. Speaking solemnly, it is a measure
defined there.

   We can flip a coin 3 times but pay attention only to what happens at first
two flips. This means that we take an evident projection \footnote{Don't
confuse it with the map $\mathbb{R}\to \mathbb{S}^1$  also denoted
by $p$.}
\begin{equation*}
p\colon\{0,1\}^3 \to \{0,1\}^2 \qquad p(a_1,a_2,a_3) = (a_1,a_2)
\end{equation*}
and pay attention only to those events --- subsets of $\{0,1\}^3$ --- which
are preimages of subsets of $\{0,1\}$ (essentially, of those events which
happened during the first two trials). Using the analogous projection
\begin{equation*}
p_1\colon \{0,1\}^3 \to \{0,1\} \qquad p(a_1,a_2,a_3)=a_1,
\end{equation*}
we can say that in the previous example with events $A,B$
$$ B = p_1^{-1} \{0\} \cap A.$$

   Idealizing the reality, we shall consider infinite sequence of a coin flips.
An elementary event is now a result of such sequence of trials; it can be
described by an infinite sequence $(a_0,a_1,a_2,\ldots)$ of symbols $0,1$. More
formally, we shall regard these sequences themselves as elementary events. It
will be convenient to us to make a slight modification of what was said and to
assume that the coin is lying before us and we see what face is above at the
moment; let $a_0$ be the number associated to this face. An elementary event
from now on is an infinite sequence $(a_0,a_1,\ldots,a_n,\ldots)$ where, once
more, $a_0$ is what we see at the very beginning (at the moment zero) and $a_n$
is the result of the $n$-th trial --- assuming that the trial is made every
second, it is what we shall see in $n$-th second. Then $\{0,1\}^{\mathbb{Z}_+}$
is the space of elementary events. Earlier we had a notion of a cylindric set.
Such sets appearing when we are fixing some coordinates, --- say, coordinates
with numbers $i_1,\ldots,i_n,$ --- correspond to the point of view when we
are interested only in what was the result not of all trials, but only
of the trials with numbers $i_1,\ldots,i_n,$. Using the evident projection
\begin{equation*}
\Omega \to \{0,1\}^n \qquad \mbox{sequence } (x_i; \ i \in \mathbb{Z}_+) \mapsto
(x_{i_1},\ldots,x_{i_n}),
\end{equation*}
we see that cylindric sets are preimages of elementary events from $\{0,1\}^n$
under this projection. (Note that $n$ can be different for
different cylindric sets.)

   Cylindric sets certainly must be considered as events (to see such and such
faces in such and such moments of time is certainly a rather elementary kind of event).
If restrictions are imposed at $n$ moments of time, the probability of the
cylindric set is $\frac{1}{2^n}$, if the coin is ``fair''. For an ``unfair''
coin the probability is $p_{a_{i_1}}\ldots p_{a_{i_n}}$, i.~e.\ if $k$ of the
numbers $a_{i_j}$ are 0 (and $n-k$ are 1), then the probability is
$p_0^kp_1^{n - k}$.
After this one can define the notion of the probability for some more
complicated subsets of $\Omega$. Essentially it is the same process which can be
used for defining the measure $\mu$ above, have not we done this differently ---
defining $\mu$ as the preimage of the standard Lebesgue measure $\mes$
under the map (\ref{pitoS1}). In any
case, for ``fair'' coin we already have a desired measure at our treatment
--- this is just $\mu$ constructed above. For an ``unfair'' coin we have
to do some work which we shall omit. By the way, in this case one can again
receive $\mu$ as the preimage of some measure on $\mathbb{S}^1$, but this
measure on $\mathbb{S}^1$ is not the well-known Lebesgue measure,
but some Lebesgue---Stieltjes measure. In many textbooks a construction of
such measure on the base of a given distribution function is described; taking
this as granted, we can easily pass to $\mu$ --- we mainly have only to
describe the distribution function which we need, and this is relatively easy.
Of course, in both cases one can avoid
going into details with $\mu$ simply because they are essentially contained
in the more well-known construction of the Lebesgue measure or of the slightly
less well-known Lebesgue-Stieltjes measure. The latter construction, which
historically was the prototype of analogous and more general constructions;
also begins from the most elementary case (``measure of an interval is its
length'') and then goes step by step to more general sets. Simplification in
our case is due to the fact that we need not imitate this construction but
can use in a formal way results of this construction carried over
on $\mathbb{S}^1$ or, what is the same, on $[0,1)$.

   And now we can finish comparing of our dynamical system with the random process of
the coin flips. A random function is a measurable function on
$\Omega$. A random process is a sequence of random functions
$\varphi_n$; $\varphi_n(x)$ is what we shall observe at the moment
$n$ provided an elementary event $x$ is realized. Denote by $\xi$
a function on $\Omega$ which is simply the projection on the
zeroth coordinate. Then the result of the $n$-th flip is
$\xi(\sigma^n(x))$. It is a sequence of numbers describing to what
of our semicircles $C_0,C_1$ comes the moving phase point (jumping
every second from $x$ to $f(x)$) at the moment $n$.

   Borel showed how the notions and facts of the measure theory\footnote{Needless
to recall that it was he who started a fruitful work towards creation of this
theory, disregarding earlier attempts which were much less satisfactory.} in order
to define in a reasonable form the notion of probability for a rather broad class of
events (subsets of $\Omega$). This allowed to study problems such that the whole
infinite sequence of trials was involved in a more essential way than before.
Borel's strong law of large numbers was the first example of this new trend, which
turned out to be fruitful. This is what we had in mind saying that Borel's impact
on the foundations of the probability theory was only one side of his work
(but, of course, these sides were closely tied).

   But at the same time Borel encountered an example of the ``chaoticity'' in the
theory of dynamical systems. This was not understood in his time --- one more
manifestation of the chaoticity in this area. The fact that there are dynamical
systems which are, so to speak, ``intrinsically chaotic'' (chaotic due to
their own dynamics, not because of exterior perturbations) and the mechanism
making them chaotic\footnote{At least the mechanism making many systems
chaotic. We do not claim that there can be no other sources of chaoticity.}
were understood much later, in 1960s.

\section{Hyperbolic automorphisms of the 2-torus}

\begin{figure}[t]
\begin{center}
\begin{tabular}{ll}
a.~\raisetonum{\includegraphics{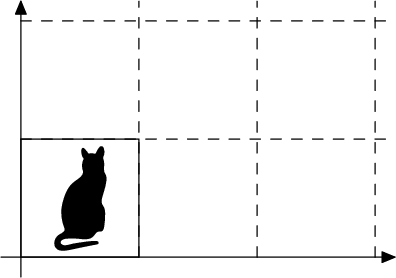}}\qquad&%
b.~\raisetonum{\includegraphics{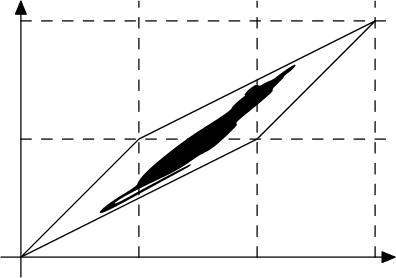}}\\
\multicolumn{2}{l}{\rule{0pt}{11pt}}\\
c.~\raisetonum{\includegraphics{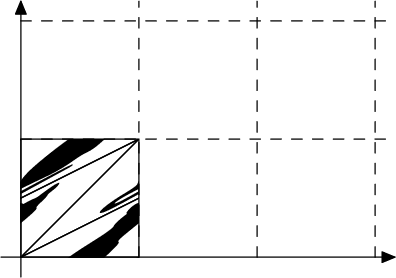}}\qquad&%
d.~\hspace{5pt}\raisetonum{\includegraphics{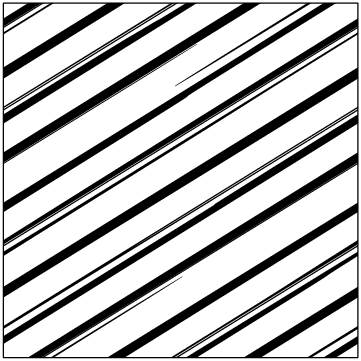}}\\
\end{tabular}
\caption{\protect\begin{tabular}[t]{@{}l@{}}
a--c. Action of $A=\protect\psmallfour2111$ on torus;\protect\\
d. Action of $A^3$ on Fig.~a, magnified.\protect\\\protect\end{tabular}}
\label{Fig2111}
\end{center}
\end{figure}
   An algebraic automorphism of the 2-torus $\mathbb{T}^2 = \mathbb{R}^2/\mathbb{Z}^2$
(the standard  projection $\mathbb{R}^2 \to \mathbb{T}^2$ will be
denoted by $p$) is defined by a matrix $A \in
\mbox{SL}(2,\mathbb{Z})$ or $A\in \mbox{GL}(2,\mathbb{Z})$.
Initially, $A$ acts on $\mathbb{R}^2$ and then this action
projects onto $\mathbb{T}^2$. Namely, $A$ defines a toric
automorphism
\begin{equation*}
\widehat{A}\colon \mathbb{T}^2 \to \mathbb{T}^2 \qquad  \widehat{A}p(x) = p(Ax),
\mbox{ i.e. } \widehat{A} (x + \mathbb{Z}^2) = Ax + \mathbb{Z}^2.
\end{equation*}
$\widehat{A}$ and $A$ are called \emph{hyperbolic} if for the eigenvalues
$\lambda,\mu$ of $A$ one has $|\lambda| > 1, \ |\mu| < 1$. Let $E^u_A$ be
the unstable eigendirection for $A$, i.e. a line $\mathbb{R} e$ in $\mathbb{R}^2$ where $Ae =
\lambda e$; later we shall also need the stable eigendirection
$E^s_A = \mathbb{R}e'$ where $Ae' = \mu e'$. Denote by $W^{u,s}_A$
the projections of $E^{u,s}_A$ to $\mathbb{T}^2$. They are dense on the torus.
Projections of the lines parallel to $E^{s,u}_A$ constitute an unstable
(expanding), resp. stable (contracting) foliation ${\cal W}^{u,s}_A$ on
$\mathbb{T}^2$; it consists of the lines obtained from $W^{u,s}_A$
under the actions of the group shifts. (We shall need ${\cal W}^{u,s}_A$
only in Parts 3 and 4.)

Figure \ref{Fig2111} is a ``standard'' illustration for the hyperbolic
automorphism of $\mathbb{T}^2$. It concerns $A=\psmallfour2111$ and presents
the action of $A$ on a figure $C$ in a fundamental square $[0,1]^2$
(Fig. \ref{Fig2111}a). Traditionally, $C$ represents a cat's silhouette,
so-called ``Arnold's cat''. On the covering plane an image of $C$ under
the action of $A$ partially leaves $[0,1]^2$ (Fig.~\ref{Fig2111}b),
so we cut it into several pieces and return them into the unit square by
shifts $(x,y)\mapsto (x+m,y+n)$ with $m,n\in\mathbb Z$ (Fig.~\ref{Fig2111}c).
Figure \ref{Fig2111}d illustrates \emph{mixing property} of this map:
for any measurable sets $X$ and $Y$ one has $\mes(\widehat{A}^n X\cap Y)\to\mes(X)\mes(Y)$
as $n\to\infty$. This means that a proportion of $Y$ occupied by $\widehat{A}^n X$
is approximately the same as the proportion of the entire torus occupied
by $X$ (equivalently, $\widehat{A}^nX$). We see that if $X=C$ and $Y$ is a quite
large rectangle then even for $n=3$ this equality holds with good precision.

   Map $\widehat{A}$ of the torus is in an evident sense expanding along ${\cal W}^u_A$
(expanding in the direction of $\mathcal{W}^u_A$), so one has the
same phenomenon of quickly increasing uncertainty as it happens
for the expanding circle map $f$ from Part 1 does. Thus it is not
surprising that the dynamical system $\{\widehat{A}^n\}$ on
$\mathbb{T}^2$ also resembles some stochastic processes.

   Many ``stochastic'' features of $\{\widehat{A}^n\}$ were revealed dealing with this system
itself. But now the most lucid way of revealing them is to use the so-called
``Markov partitions" introduced (in this case) by R.Adler and
B.Weiss\footnote{
There exists a more general version of the Markov partitions.
First step towards its elaboration was made by Ya.G.Sinay
(partially together with B.M.Gurevich), final version is due to
R.Bowen. He elaborated it for general hyperbolic sets. Subsequent
steps were to introduce (and to use) the analogous partitions
(also called ``Markov'') for several objects which are not
hyperbolic sets but which resemble them in some important aspects
--- pseudo-Anosov maps, Lorenz attractors,  some billiards ... The
works of various authors where these steps were made could be very
good, but as it concerns the general idea of the Markov partition,
essentially here we meet not so much a further development of this
general idea, but rather its adopting to a somewhat new situation.
\par
   We shall speak only about the case considered by Adler and Weiss.
It is more simple and lucid geometrically than these
generalizations and modifications. (Some exception is the
pseudo-Anosov case which is also two-dimensional and also admits
sufficiently understandable pictures. (A.Yu.Zhirov even provided
an album with such pictures --- to appear at the site of the
Steklov Inst.) But this case in more complicated in its essence
and, in our opinion, much has be done in this case before it will
become compatible
to the classical one in all respects.)}. 
They will be considered in the next part.
Here we dwell on another question. If we are interested in hyperbolic automorphisms
of $\mathbb{T}^2$, then why not to try to classify them?

   It is reasonable to consider two objects related to $\mathbb{T}^2$ as ``similar''
or ``equivalent''  if there exists a homeomorphism $\varphi: \mathbb{T}^2 \to \mathbb{T}^2$
transforming one object into another. This makes sense if we can speak about the
action of $\varphi$ on the objects considered. For the map $\widehat{A}:\mathbb{T}^2 \to \mathbb{T}^2$
it is reasonable to say that $\varphi$ transforms $\widehat{A}$ into the map
$\varphi \circ \widehat{A} \circ \varphi^{-1}$.\footnote{ As $\widehat{A}$ maps $x$ into
$\widehat{A}(x)$, it is  reasonable to say that $\varphi$ transforms $\widehat{A}$ to the map
which maps $\varphi(x)$ to $\varphi (\widehat{A}x)$,} So for the automorphisms
$\widehat{A},\widehat{B}$ of the two-torus we consider $\widehat{B}$ as ``similar'' to
$\widehat{A}$ if and only if there exists a homeomorphism $\varphi$ such that
$\widehat{B} = \varphi \circ \widehat{A} \circ \varphi^{-1}$. Then for the induced maps
\begin{equation}\label{homol}
(\widehat{B})_*, (\widehat{A})_*, \varphi_* : H_1(\mathbb{T}^2,\mathbb{Z}) \to H_1(\mathbb{T}^2,\mathbb{Z})
\end{equation}
of the one-dimensional homology group we have
$$ (\widehat{B})_* = \varphi_* \circ (\widehat{A})_* \circ \widehat{\varphi}_*^{-1} .$$
It is well known that under a suitable (and the most natural)
choice of the basis in $H_1(\mathbb{T}^2,\mathbb{Z})$ maps
(\ref{homol}) are described by matrices $A,B$ and some $C \in
\mbox{GL}(2,\mathbb{Z})$. Thus we have to deal with the usual
conjugacy of matrices $A$ and $B$. Of course now the conjugacy has
to be performed via a matrix $C$ that itself belongs to
SL($2,\mathbb{Z})$ or GL($2,\mathbb{Z})$. Conversely, if $B =
CAC^{-1}$ with $C \in \mbox{GL}(2,\mathbb{Z})$, then $\widehat{B}
= \widehat{C} \widehat{A} \widehat{C}^{-1}$. So we arrive at the
question: given hyperbolic $A$ and $B$, how to decide whether they
are conjugate in GL$(2,\mathbb{Z})$?

   If we consider a more broad conjugacy: $A \sim B$  if and only $B = CAC^{-1}$ with some
$C \in \mbox{GL}(2,\mathbb{C})$, one can find the answer in a usual course of linear algebra.
A necessary condition for such equivalence is that $A$ and $B$ have the same eigenvalues.
And if eigenvalues of a matrix are different (what is the case for our $A$ and $B$),
this condition is also sufficient. Moreover, if the eigenvalues are real (what is also the
case for our $A$, $B$), then the conjugacy can be performed via a real matrix, i.e. there
exists $C \in \mbox{GL}(2,\mathbb{R})$ such that $B = CAC^{-1}$.

   But we want to have $C \in \mbox{SL}(2,\mathbb{Z})$ or $\in \mbox{GL}(2,\mathbb{Z})$. It turns
out that this really is an additional requirement.

   This was known to Gauss. Indeed, Gauss reduced the question to the question in the
theory of binary quadratic forms. The last question was solved by him.
Now we describe this reduction.

Let $q=(A,B,C)$ be a quadratic form. For our consideration, we suppose
all coefficients of quadratic forms to be integer.
 We define its action on a vector $z=(x,y)^T$
as $q(z)=Ax^2+Bxy+Cy^2$.
Further, a \emph{discriminant} of the quadratic form $q$ is
denoted as $\mathop{\mathrm{disc}}q$ an is
equal to $B^2-4AC$. We denote by $Q(D)$ the class of all quadratic forms
with $\mathop{\mathrm{disc}}q=D$.
The group $SL_2(\mathbb Z)$ acts on $Q(D)$ by natural formula
\begin{equation*}
(g^*q)(z)=q(g^{-1}z).
\end{equation*}
On the other hands, this group acts on sets $H_\pm(t)$ of all hyperbolic
automorphisms with a given trace $t$ and a given determinant $\pm 1$
by conjugation:
\begin{equation*}
a_g\colon X\mapsto g Xg^{-1}.
\end{equation*}
Now we construct a bijection $f\colon H(t)\to Q(t^2-4)$ such that
the following diagram is commutative.
\begin{equation}\label{QFvsHA}
\arraycolsep=0pt
\begin{array}{ccc}
H_\pm(t)&{}\xrightarrow{f}{}&Q(t^2\mp4)\\
\llap{$\scriptstyle a_g$}{\downarrow}&&\llap{$\scriptstyle g^*$}{\downarrow}\\
H_\pm(t)&{}\xrightarrow{f}{}&Q(t^2\mp4)\\
\end{array}
\end{equation}
This diagram performs the desired reduction.

Now, to prove \eqref{QFvsHA}, put $f(X)(z)=\mathop{\mathrm{disc}}(\det(z,Xz))$, here $(z,Xz)$
is a $2\times 2$-matrix consisting of two columns $z$ and $Xz$.
Firstly, by direct calculation we obtain
\begin{equation*}
f\Bigl(\begin{matrix}a&b\\c&t-a\\\end{matrix}\Bigr)\,
\Bigl(\begin{matrix}x\\y\\\end{matrix}\Bigr)
=cx^2+(t-2a)xy-by^2,
\end{equation*}
so $\mathop{\mathrm{disc}}(f(X))=t^2-4\det X=t^2\mp 4$. Then, for any
form $q=(A,B,C)\in Q(t^2\mp 4)$
there exists a unique $X=\psmallfour abc{t-a}\in H_\pm(t)$ such that $f(X)=q$.
Indeed, $c=A$, $b=-C$, $a=(t-B)/2$, and to check $a$ to be integer we note
that $B^2-t^2=4AC-4$, so $B$ $t$ are of the same parity.

Finally, prove the diagram to be commutative:
\begin{multline*}
f(a_g(X))(z)=\det(z,gXg^{-1}z)=\det(g)\det(g^{-1}z, Xg^{-1}z)={}\\
{}=\det(g)\cdot f(X)(g^{-1}z)=\det(g)\cdot (g^*(f(X)))(z),
\end{multline*}
so since $\det(g)=1$, the proof is completed.

But we prefer to present an answer to our question (not only the statement
of this answer, but also the way leading to it) in terms more specific for
our framework. It seems that this rephrasing of Gauss' result and his
arguments should be well-known, but we don't know any references on this matter.

   Let $E^u_A$ be as before (the unstable eigendirection for $A$). As a line on
$\mathbb{R}^2$, it has equation $x = \kappa_A y$, with $\kappa_A$ being a
quadratic irrationality.
According to Lagrange, its continued fraction expansion
is periodic:
\begin{multline}\label{EU_Decomp}
\kappa_A = [a_0;a_1,a_2,\ldots,a_k ,\underline{a_{k + 1},\ldots,a_{k + q}},
\underline{a_{k + q + 1},\ldots,a_{k + 2q}},\ldots]={}\\
{}=[a_0;a_1,a_2,\ldots,a_k ,(a_{k + 1},\ldots,a_{k + q})]
\end{multline}
($a_{k + iq + j} = a_{k + j}$ for $i \geq 0, \ j = 1,\ldots,q$). By ``the period''
of this continued fraction we shall mean not only $q$, but also the finite sequence
of numbers $(a_{k + 1},\ldots,a_{k + q})$ up to a cyclic permutation.
The final result about the conjugacy is:

   $A$ is conjugated to $B$ via some $C \in \mbox{GL}(2,\mathbb{Z})$ if and
only if the continued fraction expansions of $\kappa_A$ and $\kappa_B$
have the same period (i.e. the same periodic part).

   Here follows a brief sketch of the proof. It is based on the following three
facts.

   a) Quadratic irrationalities $\kappa, \kappa_1$ have the same period if
and only if $\kappa_1$ can be obtained from $\kappa$ by applying to $\kappa$
some sequence of the following transformations:
\begin{equation*}
T_1(\kappa) = \kappa + 1, \quad T_2(\kappa) = \frac{1}{\kappa}, \quad
T_3(\kappa) = - \kappa
\end{equation*}
and their inverses. This easily follows from the formulas
$$ T_1([a_0;a_1,a_2,\ldots]) = [a_0 + 1;a_1,a_2,\ldots],$$
$$ T_2([a_0;a_1,a_2,a_3,\ldots]) =
\begin{cases}
[a_1;a_2,a_3,\ldots], &\text{if $a_0 > 0$,}\\
[0;a_0,a_1,a_2,\ldots],&\text{if $a_0 = 0$,}\\
\text{\rlap{(some cases for $a_0<0$),}}\\
\end{cases}$$
$$ T_3([a_0;a_1,a_2,a_3,\ldots]) =
\begin{cases}[-a_0 - 1; a_2 + 1,a_3,\ldots], &\text{if $a_1 = 1$, }\\
[-a_0 - 1;1,a_1-1,a_2,a_3,\ldots], &\text{if $a_1 \neq 1$.}
\end{cases}$$
We do not present all cases for $T_2$ due to large number of them.
This cases, where $\kappa$ is negative, can be obtained from the formula
$T_2(\kappa)=T_3(T_2(T_3(\kappa)))$. Here in the right-hand side $T_2$ is
applied to $-\kappa>0$. Note also that even in these cases $a_n$ with large
numbers shift by odd number of positions ($\pm 1$ or $\pm 3$).

  b) $\kappa_{C_i A C_i^{-1}} = T_i(\kappa_A)$, where\footnote{Here is a slightly more sophisticated point of view on the
relations between $T_i$ and $C_i$. The standard action of the nondegenerate
matrices $C = \left(\begin{array}{cc} \alpha & \beta \\ \gamma & \delta
\end{array}\right)$ on $\mathbb{R}^2$
\begin{equation*}
z = \left( \begin{array}{c} z_1 \\ z_2 \end{array} \right)
\mapsto w = \left( \begin{array}{c} w_1 \\ w_2 \end{array} \right) = Cz
\end{equation*}
defines also their action on the projective line $\mathbb{RP}^1$ considered
as the space of the straight lines passing through the origin: simply
$L \mapsto C(L)$. On
\begin{equation*}
\mathbb{RP}^1 \setminus \{\mbox{ the horisontal line } w_2 = 0 \}
\end{equation*}
we have the natural coordinate $\kappa = \kappa(L)$ that is the slope
of $L$ (so $L$ is described by the equation $z_1 = \kappa z_2$ mentioned above.
One can associate to a horizontal line the symbol
$\infty$ having in mind the usual agreements about the algebraic operations
with $\infty$.). Then for a line $L$
\begin{equation*}
\kappa (C(L)) =
\frac{\alpha \kappa(L) + \beta}{\gamma \kappa(L) + \delta}.
\end{equation*}
Denote the fractional linear transformation $\kappa \mapsto
\frac{\alpha \kappa + \beta}{\gamma \kappa + \delta}$ by $T(C)$
(we can extend it to the whole $\mathbb{RP}^1$ taking $T(C)\infty
= \frac{\alpha}{\gamma}$, but we do not need this). Then $T(C_i) =
T_i, \quad i = 1,2,3$. It remains to add that $C(E^u_A) =
E^u_{CAC^{-1}}$.}
\begin{equation*}
C_1 = \left( \begin{array}{cc} 1 & 1 \\ 0 & 1 \end{array} \right) , \quad
C_2 = \left( \begin{array}{cc} 0 & 1 \\ 1 & 0 \end{array} \right), \quad
C_3 = \left( \begin{array}{cc} -1 & 0 \\ 0 & 1 \end{array} \right).
\end{equation*}

   c) These $C_i$ are generators of GL$(2,\mathbb{Z})$.

Thus if $\kappa_A$ and $\kappa_B$ have the same period for some
$A,B \in \mbox{GL}(2,\mathbb{Z})$, then due to statement a)
$\kappa_A$ can be obtained from $\kappa_B$
by a sequence of transformations $T_i^{\pm 1}$. So $A$ is obtained from $B$ by
conjugation with a corresponding product of matrices (because of b)).

Conversely, c) implies that if $B = CAC^{-1}$ with some
$C \in \mbox{GL}(2,\mathbb{Z})$,
then $B$ can be obtained from $A$ by conjugation by some product
of $C_i^{\pm 1}$ and so $\kappa_A$ and $\kappa_B$ have the same period.

   As regards to the conjugation via $C \in \mbox{SL}(2,\mathbb{Z})$, we shall mention
only the following:

   If the period $q$ (``the length of the periodic part'') of the continued fraction
   expansion for $\kappa_A$ is odd,
and $A\sim B$ via some $C \in \mbox{GL}(2,\mathbb{Z})$, then $A \sim B$ via some
$D \in \mbox{SL}(2,\mathbb{Z})$;

   if the period is even and $A \sim B$ via some
$C \in \mbox{GL}(2,\mathbb{Z}) \setminus \mbox{SL}(2,\mathbb{Z})$, then
there is no $D \in \mbox{SL}(2,\mathbb{Z})$ conjugating $A$ and $B$.

     Both statements are simple consequences of the following ones:

     (a) if $q$ is odd, there exists a matrix $C \in \mbox{GL}(2,\mathbb{Z})$ such that
$\det C = -1$ and $A = CAC^{-1}$;

     (b) if $q$ is even and $A = CAC^{-1}$ with some
$C \in \mbox{GL}(2,\mathbb{Z})$, then $\det C = 1$.

\noindent
     Indeed, when we apply the operations $T_2$ or $T_3$ to $\kappa_A$, this leads
to a shift on one position left or right of all coefficients of the
continued fraction expansion for $\kappa_A$ with sufficiently large
number: $n$-th coefficient $a_n$ goes to the $(n + 1)$-st or $(n - 1)$-st
place. When we apply $T_1$, $a_n$ remains on the $n$-s place. Here we speak
about the ``fate'' of an individual  coefficient under the action of $T_i$
on $\kappa_A$. This needs some care, but can be justified for $a_n$ with
large $n$. On the other side, $\det C_1=1$, $\det C_2=\det C_3=1$, so
for any $C\in\mathrm{GL}(2,\mathbb Z)$
\begin{equation*}
\det C=1\iff \text{\parbox{8cm}{$C$ shifts the ``tail'' of continued fraction for $\kappa_A$
by an even number of positions.}}
\end{equation*}

So, if the period is even, then any transformation that maps $\kappa_A$ to
$\kappa_A$ should shift its ``tail'' by $qt$ ($t\in\mathbb Z$) positions that is even number.
Therefore, determinant of a corresponding matrix should be equal to $1$.

On the other hand, if this period is odd then it is not difficult to make sure
that there exists a sequence of transformations that shifts ``tail'' exactly
by $q$ positions (so, determinant of the matrix should be $-1$).

     For example, if $A =
\bigl(\begin{smallmatrix}2 & 1 \\ 1 & 1 \end{smallmatrix}\bigr)$, then
$\kappa_A = \frac{1 + \sqrt{5}}{2} = [(1)]$ and so $\kappa_A =
\frac{1}{\kappa_A - 1} = T_2 T_1^{-1}(\kappa_A)$.

Consequently,
$A = (C_2C_1^{-1})A(C_2C_1^{-1})^{-1}$ (what can be checked directly),
where $\det (C_2C_1^{-1}) = -1$.

\section[Markov partitions for hyperbolic automorphism of 2-torus]%
{Markov partitions for hyperbolic\\ automorphism of 2-torus}

   First we shall define Markov parallelograms.

  a) A Markov parallelogram in the plane (for a hyperbolic
$A \in \mbox{GL}(2,\mathbb{Z})$) is a parallelogram $\Pi$ in $\mathbb{R}^2$
having two sides parallel to $E^u_A$ (let us call
these sides ``unstable'', or ``expanding'', and denote their union
by $\partial^u\Pi$)
and two other sides parallel to $E^s_A$ (let us call these sides ``stable'',
or ``contracting'', and denote their union by $\partial^s \Pi$).

   b) A Markov parallelogram in the torus (for a hyperbolic automorphism
$\widehat{A}$) is a projection $P = p\Pi$ of some Markov
parallelogram $\Pi \subset \mathbb{R}^2$ (for the related $A$)
provided that interior $\Int\Pi$ projects
injectively.\footnote{Two opposite sides of $\Pi$ may project onto
two partially overlapping arcs.} By the ``interior'' of $P$ one
often understands the image $P^\circ = p(\Int\Pi)$ of the interior
$\Int\Pi$.\footnote{Because of what is said in the previous
footnote, $P^\circ$ may be slightly less than the true interior
$\Int P$ on torus.} Projections of the unstable (stable) sides of
$\Pi$ are called the unstable (stable) sides of $P$, their union
is denoted by $\partial^u P \quad (\partial^s P)$; so $P \setminus
P^\circ = \partial^u P \cup \partial^s P$. Unstable (stable) sides
of $P$ are arcs of the leaves of the one-dimensional foliations
${\cal W}^u_A \quad ({\cal W}^s_A)$ introduced in the beginning of
Part 2).

   A Markov partition ${\cal P} = \{P_1,\ldots,P_k\}$ (for $\widehat{A})$
is a partition of $\mathbb{T}^2$ consisting of a finite number of Markov
 parallelograms $P_i$ provided this system of parallelograms satisfies
two conditions concerning its behavior with regards to $\widehat{A}$.
These conditions are formulated below. But first we must make a warning.
Strictly speaking, ``partition'' here is not a partition in a literal sense,
i.e. a decomposition of $\mathbb{T}^2$ into a system of non-intersecting sets.
In our case this means that sides of two parallelograms can have
common points. Two unstable sides (or two stable sides) of two different
parallelograms can partially overlap, they also can have a single common
point. A stable side of one parallelogram and an unstable side of another
also can have a finite number of common points. Here is a more brief formulation of the
requirement on $P_i$: $P_i^\circ$ do not intersect each other and
$\mathbb{T}^2 \setminus (P_1^\circ \cap \ldots \cup P_k^\circ$) is a finite union of
arcs lying on leaves of ${\cal W}^{u,s}_A$. Points of this set can
be considered as exceptional ones. The set of exceptional points is negligible
in many aspects (e.g. from the measure-theoretical point of view) and at the
same time this set admits a more or less concise description and thus
can be taken into attention if necessary.

   Now we shall formulate two conditions on the behavior of $\cal P$ with
respect to $\widehat{A}$.

   \textbf{I.} Each contracting side of any $\widehat{A}P_i$ lies on a contracting
side of some $P_j$. Each expanding side of any $P_i$ lies on an expanding side
of some $\widehat{A}P_j$ (i.e. on the image of an expanding side of $P_j$).

   The same can be expressed in terms of the system of Markov parallelograms
$\Pi_i$ in $\mathbb{R}^2$ mentioned in the definition of Markov parallelograms
$P_i$ in $\mathbb{T}^2$. This version of condition~I is almost literally the same as
the version formulated in terms of $P_i$; one needs only to have in mind
that in order to get a partition of $\mathbb{R}^2$, one must take $\Pi_i + (m,n)$
with all $m,n \in \mathbb{Z}$ and $i = 1,\ldots,k$.

   Another condition can be more pictorially formulated in terms
of $\mathbb{R}^2$.

   \textbf{II.} For all $i,j = 1,\ldots,k$ only one of the intersections
$A\Pi_i \cap (\Pi_j + (m,n))$ with all $m,n \in \mathbb{Z}$ can have nonempty
interior.

\noindent
In terms of $\mathbb{T}^2$ this condition claims:

   Any nonempty $\widehat{A}P_i^\circ \cap P_j^\circ$ consists of only one connectivity
component.

   Refinements of this notion.\footnote{They concern only our case (hyperbolic
automorphisms of the 2-torus, not the Markov partitions for more general or related
objects mentioned in one of the footnotes in Part 2).}

   A) Markov partitions in the strict sense (strMp) --- the Markov partitions in
the sense as defined above.

   B) Quasi-Markov partitions (qMp). Assume we are given two different
directions in $\mathbb{R}^2$  such that the straight lines going in these
directions have irrational angular coefficients. (They are not assumed to
have any relation to any $\widehat{A}$ --- now we do not have any
$\widehat{A}$ at all.) Denote by $E^1,E^2$ the straight lines going through
$(0,0)$ in these directions. Let $W^{1,2} = p(E^{1,2})$ and let
${\cal W}^{1,2}$ be one-dimensional foliations consisting of all group shifts
of $W^{1,2}$ (i.e. obtained by projecting to $\mathbb{T}^2$ all lines parallel
to $E^{1,2}$). Replacing $E^{u,s}, W^{u,s}, {\cal W}^{u,s}$ in the part of
the  definition of the Markov parallelograms and Markov partitions preceding
I,~II by $E^{1,2},W^{1,2}, {\cal W}^{1,2}$, we get a definition of a qMp
(for the two directions given).

 Let us prove that there exists no qMp consisting of merely one element,
i.e. of one Markov parallelogram. (Later we shall see that there are qMp
consisting of two elements. Such qMp's can be considered as the simplest
ones.)

Look at any point $A$ that is a corner of this parallelogram~$P$. In a small
neighborhood of $A$ boundary of $P$ is a union of two segments, one is parallel
to $E_1$, another is parallel to $E_2$. Thus there are three possibilities:
both segments have their ends in $A$ (like in letter \textsf{L}); one pass
through $A$, another ends there (like in \textsf{T}); both pass through $A$
(like in \textsf{X}).

In the first case our parallelogram should have
an angle larger than $180^\circ$. Indeed, lift $A$ to some point $\hat A$ on
the plane, choose point close to $\hat A$ that lies in more-than-$180^\circ$
angle and then consider the lifting $\Pi$ of the parallelogram that
contains this point. Then $\Pi$ is obviously not convex.

In the second case without loss of generality we can suppose that
segment parallel to $E_2$ pass through $A$ and segment parallel to $E_1$
starts in $A$ and goes in direction we call positive. Also we arbitrarily fix
positive direction on $E_2$. Any lift $\Pi$ of the parallelogram has four
corners. Note that each corner is uniquely defined by directions
of sides (there are two possibilities for a direction of edge parallel to $E_1$
that starts at the corner and two possibilities for one parallel to $E_2$).
So we see that two corners of $\Pi$, that is, (positive $E_1$, positive $E_2$) and
(positive $E_1$, negative $E_2$) project into point $A$. Thus, difference
between their coordinates on the plane is $(i,j)\in\mathbb Z^2$.
But they share the same edge of $\Pi$, which has direction $E_1$. So,
this direction has rational slope $i/j$, that is not true.

In the third case this argumentation also works, since all corners of $\Pi$
maps to the same point $A$, hence both directions $E_{1,2}$ are rational.

   C) Pre-Markov partition (preMp). Like strMp, it is also related to some hyperbolic
automorphism $\widehat{A}$, but in its definition the condition II is omitted.

   Let ${\cal P} = \{P_1,\ldots,P_k\}$ be a strMp for $\widehat{A}$. Then $\cal P$ defines
the following coding of points of $\mathbb{T}^2$ and their trajectories.

   A point $x \in \mathbb{T}^2$ is coded by a bilaterally infinite sequence
$\{i_n; \ n \in \mathbb{Z}\}$ such that $\widehat{A}^n(x) \in P_{i_n}$ for all
$n$. Strictly speaking, this coding is univalent for the points of the set
$\bigcap_{n = - \infty}^{\infty} \widehat{A}^n (P_i^\circ \cap \ldots \cap P_k^\circ)$
which is of the ``full measure'' (its complement has the Lebesgue measure 0).
Exceptional points need some special care, like points with binary rational
cyclic coordinates in Part 1), and even more care --- now the ``good'' definition
of the coding for them involves some precautions which were absent there
(see below). But still they do not make a big harm.

     We shall describe the precautions mentioned above right now, and later
we shall explain why they are taken. The previous attempt to define the
bilateral sequence $(a_n)$ corresponding to a point $x \in \mathbb{T}^2$ is
equivalent to the following receipt:

     $x \mapsto (a_n)$ if and only if $\widehat{A}^n(x) \in P_{a_n}$ for all
$n \in \mathbb{Z}$.

\noindent
In other words,
\begin{equation}\label{attempt}
     x \mapsto (a_n) \quad \mbox{if and only if } x \in
\bigcap_{n = - N}^N \widehat{A}^{-n}(P_{a_n})
\quad \mbox{for all } N \in \mathbb{Z}_+
\end{equation}
(compare to (\ref{coding}), (\ref{defpi})). Correct definition is
\begin{equation}\label{cordef}
     x \mapsto (a_n) \quad \mbox{if and only if } x \in \clos
\Biggl(\bigcap_{n = - N}^N \widehat{A}^{-n}(P^\circ_{a_n})\Biggr)
\quad \mbox{for all } N \in \mathbb{Z}_+,
\end{equation}
where $\clos$ denotes the closure. For ``unexceptional'' points
$x \in \bigcap_{n = - \infty}^{\infty} \widehat{A}^n(P_1^\circ \cup \ldots \cup
P_k^\circ)$ this definition coincides with the previous one, but if $\widehat{A}^n
x \in \partial P_i$ for some $n,i$, then for such $x$ the new definition
is more restrictive.

   It is important that different points have different codings. Thus all what happens
in the dynamical system $(\mathbb{T}^2,\widehat{A})$ is somehow reflected in the coding.

   Codes of all points constitute some subset of
$\{1,\ldots,k\}^{\mathbb{Z}}$. It turns out that it is a so-called
Markov subset. Markov subsets themselves are defined independently
of the toric automorphisms. Here follows their definition.

   Any Markov subset corresponds to  some subset
${\cal A}\subset\{1,\ldots,k\}^2$. Pairs $(i,j) \in {\cal A}$ are called
``admissible'', other pairs --- ``forbidden''. Given $\cal A$, we define the
related Markov set $M \subset \{1,\ldots,k\}^{\mathbb{Z}}$ as a set of all doubly
(bilaterally) infinite sequences $\{i_n\}$ such that $(i_n,i_{n + 1}) \in
{\cal A}$ for all $n$. $M$ is easily seen to be a closed subset of
$\{1,\ldots,k\}^{\mathbb{Z}}$  (the latter endowed by topology similar to the
topology used in Part 1) invariant with respect to the (bilateral) topological
Bernoulli shift (also defined analogously). The pair $(M,\sigma_M)$,
where $\sigma_M$ is the restriction $\sigma_M = \sigma | M$, is called the
topological Markov shift. The probability theory and the ergodic theory supply
an extensive information about $(M, \sigma_M)$.

     For a Markov subset $M$ ``coding'' points of $\mathbb{T}^2$
a pair $(i,j)$ is admissible when $\widehat{A}(P_i^\circ) \cap
P_j^\circ \neq \emptyset$, i.e. $\Int\left(\widehat{A}(P_i) \cap
P_j \right) \neq \emptyset$. The main step of proving that $M$
actually is the Markov subset corresponding to this set of
admissible pairs is the following:
\begin{equation*}
\text{if $\widehat{A}(P_i^\circ) \cap P_j^\circ \neq \emptyset$,
$\widehat{A}(P_j^\circ) \cap P_h^\circ \neq \emptyset$, \quad then
$\widehat{A}^2 P_i^\circ \cap \widehat{A}P_j^\circ \cap P_h^\circ \neq
\emptyset$.}
\end{equation*}
If we had called ``admissible'' all those points $(i,j)$ for which
$\widehat{A}P_i \cap P_j \neq \emptyset$ (what would correspond to
(\ref{attempt})), then we would have to know that
\begin{equation}\label{wrong}
\text{if $\widehat{A}(P_i) \cap P_j\neq \emptyset$,
$\widehat{A}(P_j) \cap P_h \neq \emptyset$, \quad then
$\widehat{A}^2 P_i \cap \widehat{A}P_j \cap P_h \neq
\emptyset$.}
\end{equation}
But generally the last statement is wrong. This explains why one has to define
the coding for ``exceptional'' points according to (\ref{cordef}).

\begin{figure}[bt]
\begin{center}
\includegraphics{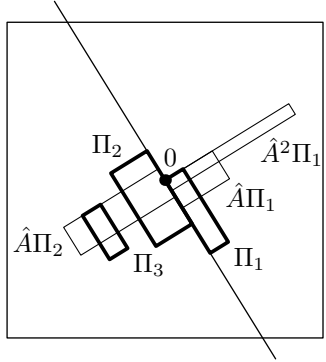}
\caption{To example showing \eqref{wrong} to be wrong.}
\label{FigForExample}
\end{center}
\end{figure}

     {\footnotesize Here is an example demonstrating that generally
(\ref{wrong}) is wrong (see Fig.~\ref{FigForExample}). For convenience we assume $\lambda$ and $\mu$ to be
positive. Denote $K = \left(- \frac{1}{2},\frac{1}{2}\right)^2$.
Clearly $(K + (m,n)) \cap K = \emptyset$, if $(m,n) \in \mathbb{Z} \setminus
\{(0,0)\}$. (The closure of $K$ is a fundamental domain.) The straight line
$E^s$ cuts $\clos K$ into two trapeziums $K'$ and $K''$ (we consider them as being
closed sets). Let Markov parallelograms $\Pi_1,\Pi_2$ be such that
$$\Pi_1 \subset K', \quad \Pi_2 \subset K'', \quad
\partial^s \Pi_1 \cap \partial^s \Pi_2 \ni 0 \ (\mbox{the origin})$$
(so that $\partial^s \Pi_i$ for both $i$ contains a small arc of $E^s$
passing through 0), and let $\Pi_1$ be so small that
$\widehat{A}^2(\Pi_1) \cup \widehat{A}(\Pi_1)\subset K'$. Finally, let $A\Pi_2$
intersect a third Markov parallelogram $\Pi_3$ lying completely in $\Int K''$.
Then
$$ \widehat{A}(P_1) \ni \widehat{A}0 = 0 \quad(\mbox{the zero of the group }
\mathbb{T}^2), \quad \widehat{A}(P_1) \cap P_2 \neq \emptyset, \quad
\widehat{A}(P_2) \cap P_3 \neq \emptyset,$$
but $\widehat{A}^2(P_1) \cap \widehat{A}(P_2) \cap P_3 = \emptyset$ and even
$\widehat{A}^2(P_1) \cap P_3 = \emptyset$, because the only ``congruent
(with respect to shifts on the elements of $\mathbb{Z}^2$) copy'' of $\Pi_3$
lying in $K$ is $\Pi_3$, which lies in $\Int K''$, while
$A^2(\Pi_1) \cap K'' = \emptyset$.

     In this argument we took as granted that there exist Markov partitions
with sufficiently small $P_i$. One can get such partition beginning with
some Markov partition and passing successfully several times from one
Markov partition to another by means of the following two operations:

     (i)~passing from a Markov partition $\{P_1,\ldots,P_k\}$ to the Markov
partition consisting of intersections $\widehat{A}P_i \cap P_j$
with nonempty interiors;

     (ii)~passing from a Markov partition $\{P_1,\ldots,P_k\}$ to the Markov
partition consisting of intersections $\widehat{A}^{-1}P_i \cap
P_j$ with nonempty interiors.}

   Originally we were interested in the dynamical system $\{\widehat{A}^n\}$
on $\mathbb{T}^2$. It turns out that the dynamical system
$\{\sigma^n_M\}$ on $M$ provides a symbolic model for the previous
system which is of the same character as the symbolic model for
$(\mathbb{S}^1,f)$ in Part 1. There exists a continuous map $\pi:
M \to \mathbb{T}^2$  such that $\pi(\mbox{the code of }x) = x$ and
$\pi \circ \sigma_M = \widehat{A} \circ \pi$. Preimage of the
Lebesgue measure on $\mathbb{T}^2$ is a measure $\mu$ on $M$
invariant with respect to $\sigma_M$. $(M,\sigma_M,\mu)$ is a
Markov process in the usual sense of the probability theory, $x =
\{x_n\} \in M$ describing the elementary event with the current
state $x_0$.

   A highly nontrivial ``purely measure theoretical'' theory of D.~Ornstein
leads to the conclusion that two Markov processes satisfying some
additional conditions which are fulfilled in our case are
isomorphic in the measure theoretical sense if (and only if ---
this was known before) they have the same entropy. Passing back to
the toric automorphisms, we can conclude that
$(\mathbb{T}^2,\widehat{A})$ and $(\mathbb{T}^2,\widehat{B})$ are
isomorphic in the measure-theoretical sense\footnote{I.~e.\ there
exists a map $\varphi: \mathbb{T}^2 \to \mathbb{T}^2$ which is an
automorphism of the measure space $(\mathbb{T}^2,\mes)$ and such
that $B = \varphi \circ A \circ \varphi^{-1}$.} if and only if
they have the same eigenvalues. (It's because the entropy in this
case is equal to $\log_2|\lambda|$ where $\lambda$ is an
eigenvalue such that $|\lambda|>1$.) Compare this with the more
complicated situation concerning the topological conjugacy of
$\widehat{A},\widehat{B}$ described in the previous part.

   Another example of the use of coding. Besides $\mu$, probability theory
provides many other measures $\nu$ which are invariant with respect to $\sigma_M$
and such that $(M,\sigma_M,\nu)$ is also a Markov process. They can be projected
to $\mathbb{T}^2$ and this supplies us with new invariant measures for $\widehat{A}$.
(While the invariance of the Lebesgue measure with respect to $\widehat{A}$
is clear, existence of other invariant measures is by no means trivial.)

   Unfortunately, the ergodic theory leads to the conclusion that usually a
strMp has to consist of rather many elements $P_i$ --- their number $k$ cannot
be less than $|\lambda|$; otherwise the diversity of motions (trajectories)
in $(\mathbb{T}^2,\widehat{A})$ cannot be reproduced in $(M,\sigma_M)$. From the
other side, any $\widehat{A}$ has a preMp consisting of two elements only.
If we shall
use this preMp for ``coding'' in the same way as it was done for a strMp,
it will turn out that two different points $x,y$ have the same coding and
the set of such $(x,y)$ is by no means ``small''. But there is a modification
of the coding process which is a remedy for this defect.

   Given a preMp $\cal P$, we define
$$ {\cal P}' = \{\mbox{closures of nonempty connected components of }
AP_i^\circ \cap P_j^\circ \}.$$
(Relations between elements of ${\cal P}$ and ${\cal P}'$ are better seen
on $\mathbb{R}^2$.) $\cal P'$ turns out to be a strMp. Thus it defines a ``good''
coding. This coding can also be seen and described in terms of $\cal P$
alone as follows. Associated with a preMp $\cal P$ there is a oriented
multigraph $\Gamma$:
\begin{itemize}
\item vertices of $\Gamma$ are parallelograms $P_i$;
\item there is an oriented edge $e$ from $P_i$ to $P_j$ if and only
if $A\Pi_i \cap (\Pi_j + (m,n))$ has nonempty interior;
\item if $\Int (A \Pi_i \cap (\Pi_j + (m,n)) \neq \emptyset$ for
several $(m,n)$, then corresponding to them there are edges going
from $P_i$ to $P_j$ (so each edge corresponds to some $P'_k \in {\cal P'}$).
\end{itemize}

   In terms of $\cal P'$, the pair $(p,q)$ is admissible if and only
if $\Int A\Pi'_p \cap (\Pi'_q + (m,n)) \neq \emptyset$ for some
$m,n \in \mathbb{Z}$.
In terms of $\Gamma$ this looks quite geometrically: the end of $e_p$ (i.e., the
edge corresponding to $P'_p$) is the beginning of $e'_q$. An infinite path in
$\Gamma$ is just a sequence of edges $\{e_{h_n}\}$  such that all pairs
$(e_{h_n},e_{h_{n + 1}})$ are admissible, i.e. that after coming to a vertex
along $e_{h_n}$, we continue our path along the edge $e_{h_{n + 1}}$.

   There exists a simple construction of the simplest preMp, i.e. those
consisting of 2 elements. Basically it is the construction of qMp consisting
of 2 elements for two directions $E^{1,2}$ with irrational angular
coefficients.

     It begins from choosing some system of data. First, it includes choosing
of an ``initial point'' $P \in \mathbb{T}^2$ (let $P = p(Q)$) and
choosing one of two lines $E^1 + Q, E^2 + Q$ which are parallel to
$E^1, E^2$ and are passing through $Q$. Let for the definiteness
$E^1 + Q$ be chosen (in the case when we choose $E^2 + Q$,
everything is going on analogously --- so to speak, $E^1$ and
$E^2$ exchange their roles). Choose one of two rays of $E^2 + Q$
beginning at $Q$ and denote it by $L$.

     Essential for the construction is an arc $I$ of $p(E^1 + Q) = W^1 + P$
which passes through $P = p(Q)$ and has endpoints $A, B$ such that

     --- $A$ is the first (after $P$) intersection of  $p(L)$ with  $I$,

     --- $B$ is the second intersection of $p(L)$ with $I$.

Let us parameterize $L$ by parameter $t$ so that (for the
definiteness) the value of $t$ corresponding to a point $z \in L$
equals to the length of the straightlinear segment $Pz$; such $z$
we denote by $z(t)$. Then our crucial condition on $A$ and $B$ is:

     $A = p(z(t_A)), \ B = p(z(t_B))$, where $T_{A,B}$ are such that
$0 < t_A < t_B$ and $p(z(t)) \notin I$ for $0 < t < t_B,\ t \neq t_A$.

\noindent Let us call this system of data --- $P, L$ and $I$ ---
the T-configuration (we think of $I$ as of the crossbar of the
letter T and of $L$ --- as of the vertical line (leg) of~T).

     {\footnotesize One needs some argument in order to prove that conditions
about the intersections of $p(L)$ with $I$ can be satisfied by means of the
proper choice of $I$. Begin with the arbitrary arc $J$ of $p(E^1 + Q) =
W^1 + P$ containing $P$ inside itself. Consider subsequent intersections
of $p(L)$ with $J$. Let them correspond to the values $t_i$ of the parameter
$t$, where $0 < t_1 < t_2 < \dots$. Note that $p(z(t_i))$ are dense on $J$.
Take
$$ i = \min \{j; \ p(t_j) \mbox{ and $p(t_{j + 1})$ lie on $J$ on the
opposite sizes of } P\}.$$ For $C,D \in J$ denote by $d(C,D)$ the
length of the arc of $J$ between points $C$ and $D$. Let
$\min\limits_{0 < j < i} d(z(t_j),P)$ be achieved at $j = h$. Then
we can take
$$ t_A = t_h, \ t_B = i + 1, \ A = z(t_A), \ B = z(t_B), \
I = \mbox{ the arc of $J$ between $A$ and $B$}.$$}

 A T-configuration defines some qMp in a natural way. Namely, let $C$ be
the next after $B$ point of the intersection of $p(L)$ and $I$ (it
is an interior point of $I$). It turns out that the arc $PC$ of
$p(E^2 + Q) = W^2 + P$ and the arc $I$ of $p(E^1 + Q) = W^1 + P$
divide $\mathbb{T}^2$ into two Markov parallelograms (for
directions of $E^{1,2}$).

To prove this we use the following idea. Move $I_j$ in the
direction $e_2$ ($e_2$ is a unit vector in $E_2$ that have the
same direction as $L$): $I_j(t)=I_j+te_2$. For small $t>0$ set
$I_j(t)\cap \pi^{-1}(D)$ contains only endpoints of $I_j(t)$. We
proceed until this holds and at some moment we have a
``catastrophe''. It is clear that ``catastrophe'' (i.~e.\ change
of the set $(I_j(t)\cap\pi^{-1}(D))-te_2$) can occur only at the
moments with $I_j(t)\cap \pi^{-1}(I)\ne\varnothing$. Such moments
are discrete (each component of $\pi^{-1}(I)$ produce at most one
such moment and only compact part, which contains finite number of
components, can contribute on a finite interval of time).
Therefore there is the first moment $t^*$ when
$\big((I_j(t)\cap\pi^{-1}(D))-te_2\big)$ changed, with two cases,
$I_j(t)\cap\pi^{-1}(D)$ is either one point or a segment. I the
first case there is no ``catastrophe'', as if for $t=t^*+\eps$ one
endpoint of $I_j(t)$ doesn't belong to $D$, then at $t=t^*$ it
coincides with $C$, and if there is a new point in
$I_j(t)\cap\pi^{-1}(D)$ for $t=t^*+\eps$ then $P$ lies in
$I(t^*)$.

In the second case we have again two possibilities: either
$I_j(t^*)\subset \pi^{-1}(I)$ or $\Int I_j(t^*)$ contains an
endpoint $z$ of $\pi^{-1}(I)$. But in the latter case $z-\eps
e_2\in D$, so $\Int I_j(t^*-\eps)\cap\pi^{-1}(D)\ne\varnothing$.
Thus, the former case takes place and $M_j=\bigcup_{0<t<t^*}\Int
I_j(t)$ is a connectivity component of $\mathbb{T}^2\setminus D$.

It remains to prove that these $M_{1,2}$ are the only connectivity
components. Consider any $z\in \mathbb{T}^2\setminus D$ and move
it in the direction $(-e_2)$ till the first intersection with $D$
at some moment $\bar t$. Then $z-\bar te_2\in\Int I_j$ for some
$j=1,2$ and therefore $z'=z-(\bar t-\eps)e_2\in M_j$. So we have a
path $\{z-\tau e_2\}_{\tau\in [0,\bar t-\eps]}$ in
$\mathbb{T}^2\setminus D$ that connects $z$ with a point in $M_j$.
Thus $z\in M_j$.

Inversely, any two-element qMp (for directions of $E^{1,2}$) can
be obtained in such way by means of a suitable T-configuration.
The proof use the same technique as the proof on non-existence of
qMp into one parallelogram.

So, we choose directions on $E^{1,2}$ in arbitrary way,
$E^{\pm,j}$ are their rays of corresponding direction started at $(0,0)$.
Also we define $W^{\pm,j}(P)=p(E^{\pm,j}+Q)$ if $P=p(Q)$.

Then we consider any point $P$ where two segments of parallelograms
boundary intersects. As before, we have three possibilities:
both have their ends here (\textsf{L});
both segments pass through $P$ (\textsf{X});
one pass through, one ends in $P$ (\textsf{T}).
Clearly, \textsf{L}-case can't take place, as one of the figures separated
by these lines has angle of more that $180^\circ$.

In \textsf{X}-case we prolong all four lines until they belongs to the
boundaries and obtain four points $P^{\pm,j}$. Note that $P^{+,1}$ belongs to
the segment of $\partial^{1}(\mathcal{P})$ that ends there and belongs to
$W^{-,1}(P^{+,1})$, and to the segment of $\partial^2(\mathcal{P})$
that passes through this point. So, near all four points $P^{\pm,j}$
the boundary has \textsf{T}-type, with directions of the
``leg'' of this \textsf{T} being different. Thus, these five
points are different. Count the corners of the parallelograms: two near
each $P^{\pm,j}$, four near $P$ (and some also may be in other points), totally
at least 12, not 8. So, this case also can't take place.

In \textsf{T}-case we can assume without loss of generality that
``leg'' of \textsf{T} belongs to $W^{+,2}(P)$. Similarly, we
obtain four different corners on the boundary: $P$, $P^{+,1}$,
$P^{-,1}$, $P^{+,2}$, and because in these points we have already
8 corners, there are no other corner on the boundary. Each segment
of the boundary has two ends, and these ends are
\textsf{T}-points, which are different for different segments. So,
boundary consists of two segments: $I=P^{-,1}P^{+,1}$ on
$E^1$-direction and $PP^{+,2}$ in $E^2$-direction. So, points
$P^{\pm,1}$ lies on $PP^{+,2}$. It is clear that $P$,
$L=W^{+,2}(P)$ $I$ comprise T-construction that produces given
qMp.

     If we are given a hyperbolic automorphism $\widehat{A}$ of $\mathbb{T}^2$,
then this construction with $E^1 = E^s, \ E^2 = E^u$ or $E^1 = E^u, \
E^2 = E^s$ gives a preMp for $\widehat{A}$, provided that $P$ is a fixpoint
for $\widehat{A}$.

\section{Classification of the simplest preMp}

     Besides the conjugating of toric automorphisms by means of toric automorphisms,
we shall consider their conjugating by means of affine diffeomorphisms of
$\mathbb{T}^2$, i.e. by means of maps
$$ z \mapsto \widehat{C}(z) = \widehat{B}z + g,$$
where $\widehat{B}$ are toric automorphisms and $g \in
\mathbb{T}^2$. In other words, $\widehat{C}$ is obtained by
projecting to $\mathbb{T}^2$ an affine map of the plane --- a map
$$ z \mapsto C(z) = Bz + b \qquad \mbox{with } B \in \mbox{GL}(2,\mathbb{Z})
\ \mbox{and } b \in p^{-1}(g).$$ We shall need only the case when
the result of the conjugating of a toric automorphism
$\widehat{A}$ by means of $\widehat{C}$ is a toric automorphism
again (actually we shall demand even more). It is easy to see that
this is the case if and only if $B^{-1}b$ is a fixpoint of $A$.

   If $\widehat{C}$ acts on the objects $O$
from some class of objects $\{O\}$, then it is natural to say that the pair
$$ (\widehat{A},\mbox{ an object $O$ somehow related to } \widehat{A})$$
is equivalent to $(\widehat{C}\widehat{A}\widehat{C}^{-1}, \widehat{C}(O))$
(provided it is true that $\widehat{C}(O)$ is related to
$\widehat{C}\widehat{A}\widehat{C}^{-1}$ in the same way as $O$ is related to
$\widehat{A}$).

   If ${\cal P} = \{P_i\}$ is a preMp for $\widehat{A}$, then
$\widehat{C}{\cal P} = \{\widehat{C}P_i\}$ is a preMp for
$\widehat{C}\widehat{A}\widehat{C}^{-1}$:

   if sides of $\Pi_i$ are parallel to $E^{u,s}$, then sides of
$CAC^{-1}(C\Pi_i)$ are parallel to $E^u_{CABC{-1}} = BE^u_A, \
E^s_{CAC^{-1}} = BE^s_A$;

   if $\widehat{A} P_i \cap P_j$ are ``good'', then
$$\widehat{C}\widehat{A}\widehat{C}^{-1}(\widehat{C}P_i)\cap
\widehat{C}P_j = \widehat{C}(\widehat{A}P_i \cap P_j)$$ are also ``good''.

From this point we impose an additional condition on preMp. Since
a contracting segment of its boundary maps into itself, there is a
fixed point on it (as segment is compact). Due to the same reason
for inverse transform the expanding segment also has a fixed
point. In our examples these two fixpoints are the same one placed
in one of the four joint points (``vertexes'') of contracting and
expanding segments, i. e. the following condition holds:

\textbf{III.} There is a fixpoint that belongs to an intersection
of stable and unstable segments.

We call such preMp's to be ``of vertex type''. There are also
preMp's without this condition with different fixpoints on
expanding and contracting segments, they are called to be ``of
edge type''. Vertex-type preMp's appears to be a source for
description of all preMp's, this will be discussed at the end of
this Part.

So, from now on until near the end of this Part, we will consider
only vertex preMp's without any special mention.

   Let $\widehat{C}\widehat{A}\widehat{C}^{-1} = \widehat{A}$ (what means
that $BAB^{-1} = A$, i.e. $B$ commutes with $A$, and $b$ is a fixpoint of
$A$). In this case we consider a preMp $\cal P$ and a preMp
$\widehat{C}{\cal P}$ as equivalent ones. Question: What is the number of
the equivalence classes of the simplest preMp for $\widehat{A}$?
Answer is given by the following theorem.

\begin{theorem}\label{NumOfClasses}In terms of \eqref{EU_Decomp} (see Part 2), there are
\begin{equation*}
2(a_{k + 1} + \ldots + a_{k + q}) = 2(\mbox{sum of the $a_i$ in the period}).
\end{equation*}
classes of (vertex) preMp's, $2q$ (twice the length of the period) of them
are of the ``island'' type, others are of the ``parquet'' type.
\end{theorem}

\begin{figure}[t]
\begin{center}
\begin{tabular}{ll}
a.~\raisetonum{\includegraphics{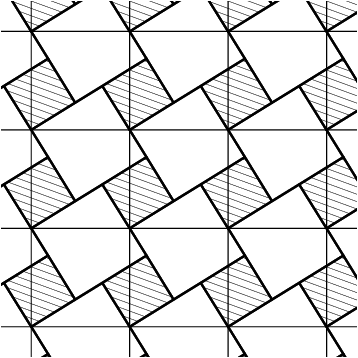}}\qquad&%
b.~\raisetonum{\includegraphics{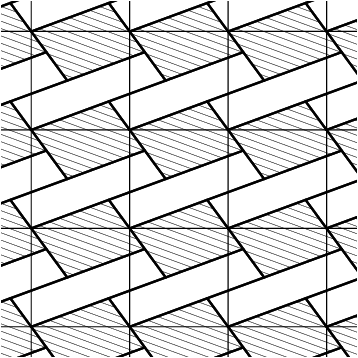}}\\
\end{tabular}
\caption{``island'' (a) and ``parquet'' (b) types of preMP's.}
\label{FigTwoTypes}
\end{center}
\end{figure}

Two types mentioned in the theorem differs by topological properties of their
lifting to the plane. For ``island'' type
there are parallelograms which are bigger
``in all directions'' (let it be $\Pi_1 + (m,n)$) and they
constitute a connected set (``ocean'' $\bigcup_{m,n} (\Pi_1 + (m,n))$);
a union of other parallelograms $\bigcup_{m,n} (\Pi_2 + (m,n))$
is disconnected and its connected components
are these  $\Pi_2 + (m,n)$ (``islands''). (See Figure~\ref{FigTwoTypes}a.)

For ``parquet'' type preMp both sets $\bigcup_{m,n} (\Pi_1 + (m,n))$ and
$\bigcup_{m,n} (\Pi_2 + (m,n))$ have infinitely
many connected components each consisting of infinitely many parallelograms;
each component resembles a stripe. (See Figure~\ref{FigTwoTypes}b.)

   In the textbooks one can meet only the island type preMp. This is because
the standard example there is
$A = \left( \begin{array}{cc}2 & 1 \\ 1 & 1 \end{array} \right)$. In this
case $\kappa_A = \frac{1 + \sqrt{5}}{2}$ (the golden mean). Its continued
fraction expansion is $[(1)]=[1;1,1,\dots]$. So there are 2 simplest preMp's
of the island type and no simplest preMp's of
the parquet type.\footnote{Note that $\psmallfour2111={\psmallfour1110}^2$,
so each of equivalence classes with respect to centralizer is split into
two equivalence classes with respect to the group $\{\pm A^n\mid n\in\Z\}$.}

   As far as we know, first picture with preMp of the parquet type was published
by E.~Rykken. But, as far as we understand, she did not discuss
when such preMp's can appear.

Now we get an outline of a proof of this result.

First, we can consider only partitions with fixpoint from condition III being
an origin~$O$. (For a shift of the torus to any vector from any fixpoint to
another one commutes with the transform.)

Further, at a small neighborhood of $O$ boundaries forms two
segments, one passes through $O$, another has its end there. So we
have four broad classes of preMp's distinguished by a direction of
the latter segment ($e_u$, $e_s$, $-e_u$, $-e_s$). But all preMp's
from the last two classes are equivalent to preMp's for the first
two of them by an automorphism $-\mathrm{id}$.

So, let us consider one of the first two classes, say $e_u$-class. We are going
to prove that there are $S$ equivalence classes, $L$ of which are of ``island'' type,
in this broad class.

\begin{lemma}\label{lem41}All preMp's from the broad class form a double infinite
sequence
\begin{equation}\label{MP_Seq}
\dots, P_{-1}, P_0, P_1, P_2,\dots
\end{equation}
such that for their stable and unstable boundary segments $I_k^{u,s}$ following
statement holds:
\begin{equation*}
I^u_k\subset I^u_{k+1},\qquad I^s_k\supset I^s_{k+1}.
\end{equation*}
\end{lemma}

\begin{proof}Let $x(t)$ be a solution of $\dot x=e_u$ with $x(0)=0$ (so $x(t)$ is
a point moving along $W^u(O)$ with a constant velocity). Denote by $(t_n)$
a sequence of all instants of time $t>0$ when $x(t)\in I$. Here $I\subset W^s(O)$
is a starting segment in T-construction. In this terms we can easily describe
a T-construction applied to any $J\subset I$. Indeed, a points $A_J$ and $B_J$
can be described as the points $x(t_{n_A})\in J$ and $x(t_{n_B})\in J$, $n_A<n_B$
with a following properties:
\begin{subequations}\label{CondNab}
\begin{align}
{}\qquad{}&\text{There are no $n<n_B$ such that}\notag\\
\label{CondNab1}
&\text{\qquad\quad $x(t_n)$ lies on $I$ between
$x(t_{n_A})$ and $x(t_{n_B})$.}\\
&\text{There are no $m<n_B$ such that}\notag\\
\label{CondNab2}
&\text{\qquad\quad$x(t_m)\in J$ and $x(t_{n_B})$
lies on $I$ between $x(t_m)$ and $O$.}
\end{align}
\end{subequations}
So, if $P_{(1)}$ and $P_{(2)}$ are two preMp's and $I=I^s_{(1)}\cup I^s_{(2)}$
we can apply this to $J=I^s_{(1)}$ and $J=I^s_{(2)}$. Without loss of generality
$n_{B_1}<n_{B_2}$ (hence $I^s_{(1)}\subset I^s_{(2)}$). So $x(t_{n_{A_1}})$
and $x(t_{n_{B_1}})$ can't lie between
$A_2$ and $B_2$, whence $I^u_{(1)}\supset I^s_{(2)}$. Thus an order
\begin{equation*}
P_{(1)}\succ P_{(2)}\iff I^u_{(1)}\subset I^u_{(2)}
\end{equation*}
is linear. Moreover, each preMp $P_{(1)}\succ P$ corresponds to
some number $n_B$ from conditions \eqref{CondNab}. So, any ``right
tail'' $(\{P'\mid P'\succ P_0\},\succ)$ is isomorphic as ordered
set to $(\N,>)$. Then the entire set of preMp's is isomorphic
either to $(\N,>)$ or to $(\Z,>)$. The former case is eliminated
due to absence of an initial element in the order: for quite long
(in both directions) initial segment $I$ the segment $AB$
corresponding to it is arbitrary long (due to density of
$W^s(O)$).
\end{proof}

\begin{lemma}$\widehat A$ (or $-\widehat A$ if $\lambda_u<0$) acts on the sequence
\eqref{MP_Seq} as a shift: $\widehat A(P_k)=P_{k+s}$.
\end{lemma}

\begin{proof}$\widehat A$ conserves the order $\succ$. Shifts are the only
automorphisms of the ordered set $(\Z,>)$.\end{proof}

For further we need to consider a structure of a centralizer of $A$ i.~e.\
a group $C(A)=\{B\in GL_2(\Z)\mid AB=BA\}$.

\begin{figure}[tp]
\begin{center}
a.~\raisetonum{\includegraphics{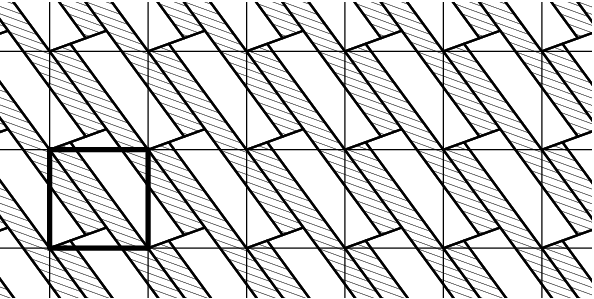}}\\[3mm]
b.~\raisetonum{\includegraphics{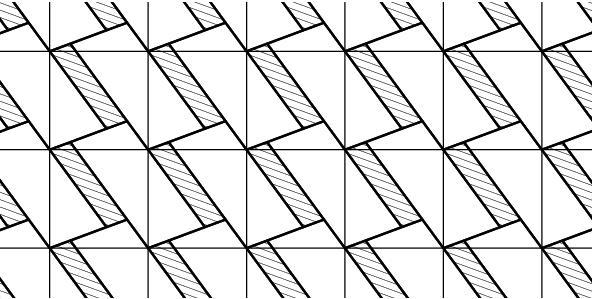}}\\[3mm]
c.~\raisetonum{\includegraphics{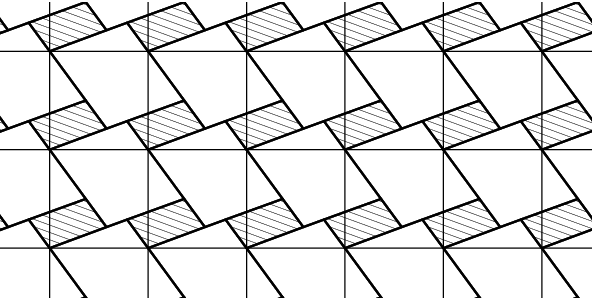}}\\[3mm]
d.~\raisetonum{\includegraphics{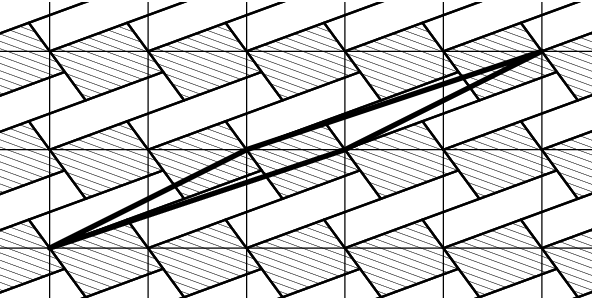}}\\
\caption{Four consecutive preMp's for $A=\protect\psmallfour3211$.}
\label{FigFourPreMps}
\end{center}
\end{figure}

\begin{lemma}\label{CentrStruct}Suppose that $A\in GL_2(\Z)$ is a hyperbolic matrix.
Then there exists $B\in GL_2(\Z)$ such that $C(A)=\{\pm B^n\mid n\in\Z\}$.
\end{lemma}

\begin{proof}There exists a matrix $D\in GL_2(\R)$ such that
$\tilde A=D^{-1}AD=\psmallfour{\lambda_u}00{\lambda_s}$.
Then  $X=D^{-1}(C(A))D$ is a subset of a centralizer of $\tilde A$
in $GL_2(\R)$, which is equal to $\bigl\{\psmallfour\lambda00\mu \mid \lambda,\mu\in \R^*\bigr\}$.
Since a conjugacy $M\mapsto D^{-1}MD$ is a homeomorphism of $GL_2(\R)$,
$X$ is a discrete set. Moreover, as $\det M=\pm 1$ for $M\in GL_2(\Z)$
this set is a subset of
$Y=\bigl\{\psmallfour\lambda00\mu \mid \lambda\mu=\pm 1\bigr\}$.
Therefore, a projection $\pi\colon\psmallfour\lambda00\mu \mapsto\lambda$
is 2:1-map, so $\pi(X)$ is a discrete subgroup of $\R^*$.

So we have two possibilities: $\pi(X)=\{\alpha^n\mid n\in\Z\}$
or $\pi(X)=\{\pm\alpha^n\mid n\in\Z\}$. The former can't take place because
$\psmallfour{-1}00{-1}\in X$. Lifting of the latter to $Y$ yields either
$X=\bigl\{\pm\psmallfour{\alpha^n}00{\beta^n}\mid n\in\Z\bigr\}$
or $X=\bigl\{\psmallfour{\pm\alpha^n}00{\pm\beta^n}\mid n\in\Z\bigr\}$
(signs are independent).

Suppose the latter case takes place. Then $F=\psmallfour100{-1}\in X$.
Therefore, $DFD^{-1}\in G(A)\subset GL_2(\Z)$.
But $DFD^{-1}$ has $e_u$ as an eigenvector with eigenvalue equal to $1$.
This means that the ratio of its coordinates should be rational,
so we have a contradiction.

Thus, $X=\bigl\{\pm\psmallfour\alpha00\beta ^n\mid n\in\Z\bigr\}$ and
the statement of the lemma is true for $B=D\psmallfour\alpha00\beta D^{-1}$.
\end{proof}

Matrix $B$ from the statement of the previous lemma can be easily
described in terms of continued fractions.

\begin{lemma}Let $e_u=(\omega,1)$, $\omega=[b_0,\dots,b_{n-1},(a_1,\dots,a_L)]$.%
\footnote{We also define $b_k$ for $k\ge n$ as follows: $\omega=[b_0,\dots,b_{n-1},b_n,b_{n+1},\dots]$.}
Then $B$ from Lemma~\ref{CentrStruct} can be chosen equal to $CDC^{-1}$,
where
\begin{equation*}
C=C_1^{b_0}T_2C_1^{b_1}C_2C_1^{b_2}C_2\dots C_1^{b_{n-1}}C_2,\qquad
D=C_1^{a_1}C_2C_1^{a_2}C_2\dots C_1^{a_L}C_2.
\end{equation*}
(Matrices $C_{1,2}$, which correspond to elementary
operations $T_1(\omega)=\omega+1$ and $T_2(\omega)=1/\omega$, were defined
in Part 2.)
\end{lemma}

\begin{proof}Denote $CDC^{-1}$ by $B'$. We can see that $e_u$ is an eigenvector
of $B'$, so $e_s$ is also an eigenvector (since they are algebraically conjugated,
as well as their eigenvalues), so $B'$ commutes with $A$.

Each matrix in $C(A)$ acts on continued fraction of $\omega$ as a shift, and
the map $d\colon C(A)\to L\Z$ that maps a matrix to the magnitude of
the corresponding shift is a group homomorphism. As $B'$ maps to $L$, $d$
should be an epimorphism. Thus $B$ should maps to $L$ or to $-L$. Then
$d^{-1}(L)=\{\pm B\}$ in the former case and $d^{-1}(L)=\{\pm B^{-1}\}$ in
the latter one. In all cases $B'=\pm B^{\pm 1}$, so $C(A)=\{\pm B'^n\mid n\in\Z\}$.
\end{proof}

Now we pass to a central point of the proof: an interrelation between
the continued fraction of $\omega$ and preMp's.

\begin{lemma}1. Let a starting segment $I$ of T-construction be
sufficiently short. Then all preMp's with $I^s\subset I$ can be
described as follows. If $A', B'$ are lifts of $A$ and $B$ that belongs to
$W^u(0,0)$ then $A'$ (correspondingly, $B'$) lies on lifts of $I\subset W^s(O)$
that consist $(p_k,q_k)$ \textup(corr., $(lp_k+p_{k-1},lq_k+q_{k-1}))$, where
$1\le l\le b_{k+1}$ and $p_n/q_n=[b_0,\dots,b_n]$ is $n$-th convergent for $\omega$.
Conversely, each such pair of points for sufficiently large $k$ corresponds
to some preMp.

2. $k$'s and $l$'s for preMp's will be arranged in \eqref{MP_Seq} as follows:
\begin{equation*}
\dots, (k-1,b_k), (k,1), (k,2),\dots, (k,b_{k+1}), (k+1,1), \dots, (k+1,b_{k+2}),\dots
\end{equation*}

3. preMp is of ``island'' type iff it corresponds to $(k,l)$ with $l=b_{k+1}$.

4. $B'$ acts on a sequence \eqref{MP_Seq} as a shift to $S=a_1+\dots+a_L$ positions.
\end{lemma}

Figure~\ref{FigFourPreMps} illustrates this lemma. There are four
consequent members of sequence~\eqref{MP_Seq} for $A=\psmallfour3211$
(here $\kappa=[0,(2,1)]$). One can see that Fig.~\ref{FigFourPreMps}d presents
the image of preMp from Fig.~\ref{FigFourPreMps}a under $A$ (the bold
parallelogram is an image of the unit square). Thus $A$ shifts
sequence~\eqref{MP_Seq} to three positions, and ``islands'' and ``parquets''
form a sequence $(P,I,I)=\ldots,P,I,I,P,I,I,P,I,I,\ldots$ as it follows from
the statements of the lemma.

Note also that the last two statements imply that there are exactly
$S$~equivalence classes (we recall that now only $e_u$-type preMp's are considered),
$L$~of~them comprises of ``island''-type preMp's. Their link to $e_s$-type
preMp's will finish the proof by Lemma~\ref{TypesEquiv} below.

\begin{figure}[t]
\begin{center}
a.~\raisetonum{\includegraphics{fig5_a.ps}}
\par
\vspace{5mm}
b.~\raisetonum{\includegraphics{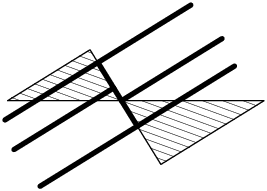}}\quad
c.~\raisetonum{\includegraphics{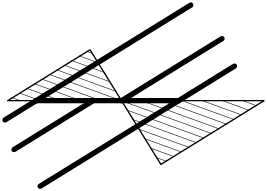}}\\
\parbox{0.7\textwidth}{\caption{The ``butterfly'' (a) and transformation of $(e_u,I)$-qMp (b) into
$(e_u,J)$-qMp (c).}
\label{FigButterfly}}
\end{center}
\end{figure}

\begin{proof}Let $I$ be so small that different ``butterflies'' on the plane
don't intersect. Here ``butterfly'' is defined as a union of two triangles
(with their interior), the boundary of each consists a connected component
of $I\setminus\{O\}$, a horizontal segment passing through $O$ and segment parallel
to $e_u$ (see Figure \ref{FigButterfly}a).

Thus there is a 1:1-correspondence between qMp's generated by T-con\-struc\-tion
for $I$ and those generated by T-construction for a horizontal segment of the
``butterfly''. (It is denoted by $J$.) This correspondence is shown on
Figures~\ref{FigButterfly}\mbox{b--c}.
It is well-defined since all transformations are inside the ``butterfly'',
which is injectively mapped into plane. Note also that the relation between
$OA$ and $OB$ is the same as one between $OA'$ and $OB'$, this will be useful
to find a type of the partition.

By the same reasoning as in the proof of Lemma 1, one can obtain that
points $A$ and $B$ for any preMp with $I^s\subset I$ are $x(t_{n_{A,B}})$
that satisfy condition \eqref{CondNab1}. As $e_u=(\omega,1)$
and $J$ belongs to an $x$-axis, all $t_n$ are integers. So, this condition
can be reformulated as such: $x$-coordinates of $A$ and $B$ are equal to
$q_{A,B}\omega-p_{A,B}$ (with $q_A<q_B$) such that
\begin{subequations}\label{CondPQ}
\begin{gather}
0\in[q_{A}\omega-p_{A},q_{B}\omega-p_{B}];\\
\text{there are no $(p',q')$ with $q'<q$ such that }q'\omega-p'\in[q_{A}\omega-p_{A},q_{B}\omega-p_{B}].
\end{gather}
\end{subequations}

Consequently, both $(p_A,q_A)$ and $(p_B,q_B)$ satisfies a following
condition:
\begin{equation}\label{LR2type}
\text{there are no $(p',q')$ such that $0<q'<q$ and $q'\omega-p'\in [0,q\omega-p]$.}
\end{equation}
Such pairs $(p,q)$ (or, more commonly, fractions $p/q$) are called \emph{one-sided
best approximations to $\omega$ of second type}. Similarly, pairs $(p,q)$
satisfying a condition
\begin{equation}\label{2type}
\text{there are no $(p',q')$ such that $0<q'<q$ and $|q'\omega-p'|<|q\omega-p|$,}
\end{equation}
are called \emph{(two-sided)
best approximations to $\omega$ of second type}.

We state a theorem from number theory describing them.

\begin{theorem}\label{Approx}1. If $\omega=[b_0,b_1,\dots,]$ then one-sided approximations
are
$ p/q=[b_0,\dots b_{k-1},l]$, where $1\le l\le b_k$. They are arranged as
\begin{equation}\label{LR2formula}
\underbrace{[1], [2], \dots, [b_0]}_{\text{from below}},
\underbrace{[b_0,1], \dots, [b_0,b_1]}_{\text{from above}},
\underbrace{[b_0,b_1,1], \dots, [b_0,b_1,b_2]}_{\text{from below}},\dots
\end{equation}
with denominators growing in the sequence.\\
2. Two-sided approximations are only the following ones:
\begin{equation}\label{TS2formula}
[b_0], [b_0,b_1], [b_0,b_1,b_2], \dots, [b_0,b_1,\dots, b_n], \dots
\end{equation}
\end{theorem}

This theorem seems to be well-known and can be proved in the way similar
to the classical theorem on two-sided approximations (see, e.g., [Kh]).

Thus, $p_A/q_A$ and $p_B/q_B$ are fractions from \eqref{LR2formula}.
However condition \eqref{CondPQ} is stronger. Obviously it can be expressed
as such: there is no approximations from the same side as $p_A/q_A$ between
$p_A/q_A$ and $p_B/q_B$ in sequence \eqref{LR2formula}.

Consequently,
\begin{equation}\label{CFofAB}
p_A/q_A=[b_0,b_1,\dots,b_k],\qquad p_B/q_B=[b_0,\dots,b_k,l],
\end{equation}
where $1\le l\le b_{k+1}$. This proves the first two statements of the lemma.
(Actually it remains to prove that
\begin{equation}\label{Intermed}
[b_0,\dots,b_k]=\frac{p_k l+p_{k-1}}{q_k l+q_{k-1}}.
\end{equation}
This can be done by induction over $k$.)

Third statement is also simple. A parallelogram with its base on
the segment $OA$ has height $q_B$ and one with base on $OB$ is of
height $q_A$. Thus if $OA'>OB'$ this preMp is of ``island'' type
and otherwise it is of ``parquet'' type. (Recall that when we
return back to $AB$ segment a type of the partition remains the
same.) Statement 2 of Theorem~\ref{Approx} implies that the former
case takes place only if $l=b_{k+1}$ in \eqref{CFofAB}.

Fourth statement of the lemma obviously follows from a fact that
$B'$ maps $(p_k,q_k)$ to $(p_{k+L},q_{k+L})$.
\end{proof}

To finish the proof of Theorem~\ref{NumOfClasses} it remains to proof
that the number of $e_u$-type preMp classes are equal to the number
of those of $e_s$-type. If trivially follows from the next (and the final one)
lemma.

\begin{lemma}\label{TypesEquiv}PreMp's of $e_u$-type and of $e_s$-type can be bijectively
corresponded in such a way that any preMp can be mapped to its
correspondent by a shift on the torus.\end{lemma}

\begin{proof}Each preMp has 4 joint points on its boundary (one of each type).
So we should just shift it to place the required joint point to the origin.
The result will be preMp, so we define two mutually inverse maps (one from
$e_u$-preMp's to $e_s$-preMp's, another is reverse). So there is
a 1:1-correspondence.\end{proof}

Now we will shortly discuss a preMp's with two different fixpoints on
the boundary. They really appears at least for some automorphisms.
For example, let us consider a standard
$\bigl(\begin{smallmatrix}2&1\\1&1\end{smallmatrix}\bigr)$-automorphism $A$ and
its large degree $B=A^N$. Then $B$ has quite many fixpoints,
which are quite densely placed on torus. Now get any preMp (for $A$)
of $e_u$-type and shift it to vectors $-\eps e_u$. If fixpoints are
quite densely placed on torus, for a rather small $\eps$ the stable segment
of shifted preMp will pass through a fixpoint. On the other hand, as this shift
is quite small, the origin will retain on an unstable segment.

Similarly to Lemma~\ref{TypesEquiv} it can be proved that any preMp
(with an arbitrary position of its fixpoints) can be obtained from, say,
some $e_u$-type preMp by some shift. The number of preMp's obtained from one
can be found algorithmically as the number of points of a lattice in
a parallelogram. Indeed, if $P$ is a joint point of the $e_u$-type, and
$U=P+xe_u$ and $S=P+ye_s$ are fixpoints then $xe_u-ye_s$ belongs to a lattice
of all fixpoints. Thus we have a parallelogram of points of the form $xe_u-ye_s$
(as $x$ and $y$ are restricted to some segments) and each point
of fixpoints lattice corresponds to a preMp.

-----------------------------------------------------------------------------

   Authors whose  results and/or ideas are used or mentioned here:

   J.~Bernoulli, J.~Lagrange, C.~Gauss, H.~Poincar\'e. A.A.~Markov (senior),
G.~Frobenius, J.~Hadamard, E.~Borel, D.~Ornstein, R.~Adler,
B.~Weiss, A.Yu.~Zhirov, E.~Rykken, D.V.~Anosov, A.V.~Klimenko,
G.~Kolutsky.

The author invisibly presented here: S.~Smale. (He was of the
major influence in the hyperbolic theory during 60s and the
beginning of 70s, but none of his works or ideas are used here
explicitly.)

\section*{References}

\parindent=0pt

[ATW] R. Adler, C. Tresser, P. A. Worfolk, Topological conjugacy of linear
endomorphisms of the 2-torus, Trans. Amer. Math. Soc. 349 (1997) no. 4,
1633-1652.

[K-H] A. Katok, B. Hasselblatt. Introduction to the Modern Theory of Dynamical
Systems. Cambridge University Press, 1995. (Encyclopedia of Mathematics and Its Applications,
Vol.~54)

[Kh] A. Khinchin, Continued Fractions. Mineola, N.Y.: Dover Publications, 1997.

\end{document}